\documentclass{amsart}

\newtheorem{theorem}[equation]{Theorem}
\newtheorem{lemma}[equation]{Lemma}
\newtheorem{corollary}[equation]{Corollary}
\newtheorem{proposition}[equation]{Proposition}

\numberwithin{equation}{section}

\newcommand{\updot}{\textstyle\cdot}

\usepackage{amscd}

\begin{document}

\title{Exponential sums on ${\bf A}^n$, II}
\author{Alan Adolphson}
\address{Department of Mathematics\\
Oklahoma State University\\
Stillwater, Oklahoma 74078}
\email{adolphs@math.okstate.edu}
\thanks{The first author was supported in part by NSA Grant
  \#MDA904-97-1-0068} 
\author{Steven Sperber}
\address{School of Mathematics\\
University of Minnesota\\
Minneapolis, Minnesota 55455}
\email{sperber@math.umn.edu}
\date{}
\begin{abstract}
We prove a vanishing theorem for the $p$-adic cohomology of 
exponential sums on ${\bf A}^n$.  In particular, we obtain new classes
of exponential sums on ${\bf A}^n$ that have a single nonvanishing
$p$-adic cohomology group.  The dimension of this cohomology group
equals a sum of Milnor numbers.
\end{abstract}
\maketitle

\section{Introduction}

Let $p$ be a prime number, $q=p^a$, and ${\bf F}_q$ the finite field
of $q$ elements.  Associated to a polynomial $f\in{\bf
  F}_q[x_1,\ldots,x_n]$ and a nontrivial additive character $\Psi:{\bf
  F}_q\rightarrow{\bf C}^{\times}$ are exponential sums
\begin{equation}
S({\bf A}^n({\bf F}_{q^i}),f)=\sum_{x_1,\ldots,x_n\in{\bf F}_{q^i}}
\Psi({\rm Trace}_{{\bf F}_{q^i}/{\bf F}_q}f(x_1,\ldots,x_n))
\end{equation}
and an $L$-function
\begin{equation} L({\bf A}^n,f;t)=\exp\biggl(\sum_{i=1}^{\infty}
  S({\bf A}^n({\bf F}_{q^i}),f)\frac{t^i}{i}\biggr).
\end{equation}
One of the basic results on exponential sums is the following
theorem of Deligne\cite[Th\'{e}or\`{e}me 8.4]{DE1}.  Let $\delta=\deg f$ and
write
\begin{equation}
f=f^{(\delta)}+f^{(\delta-1)}+\cdots+f^{(0)},
\end{equation}
where $f^{(j)}$ is homogeneous of degree $j$.
\begin{theorem}
Suppose $(p,\delta)=1$ and $f^{(\delta)}=0$ defines a smooth
hypersurface in~${\bf P}^{n-1}$.  Then $L({\bf A}^n,f;t)^{(-1)^{n+1}}$
is a polynomial of degree $(\delta-1)^n$, all of whose reciprocal
roots have absolute value $q^{n/2}$.
\end{theorem}

For exponential sums on ${\bf A}^n$, several generalizations of
Deligne's result have been proved (\cite{AS1,AS2,DL,GA}).  In all
these theorems, the hypothesis implies that $f$, regarded as a
function from ${\bf A}^n$ to ${\bf A}^1$, has only finitely many
critical points and the degree of the polynomial $L({\bf
  A}^n,f;t)^{(-1)^{n+1}}$ equals the sum of the Milnor numbers of
those critical points.  We examine these critical points more closely.

We write ${\bf F}_q[x]$ for ${\bf F}_q[x_1,\ldots,x_n]$ and consider
the complex $(\Omega^{\updot}_{{\bf F}_q[x]/{\bf F}_q},\phi_f)$, 
where $\Omega^k_{{\bf F}_q[x]/{\bf F}_q}$ denotes the module of
differential $k$-forms of ${\bf F}_q[x_1,\ldots,x_n]$ over~${\bf F}_q$
and $\phi_f:\Omega^k_{{\bf F}_q[x]/{\bf F}_q}\rightarrow\Omega^{k+1}_{{\bf
    F}_q[x]/{\bf F}_q}$ is defined by
\[ \phi_f(\omega)=df\wedge\omega, \]
where $d:\Omega^k_{{\bf F}_q[x]/{\bf
    F}_q}\rightarrow\Omega^{k+1}_{{\bf F}_q[x]/{\bf F}_q}$ is the
exterior derivative.  The map $f:{\bf A}^n\rightarrow{\bf A}^1$ has
only isolated critical points if and only if
\begin{equation}
H^i(\Omega^{\updot}_{{\bf F}_q[x]/{\bf F}_q},\phi_f)=0\qquad\text{for
  $i\neq n$},  
\end{equation}
which implies that $\dim_{{\bf F}_q} H^n(\Omega^{\updot}_{{\bf
    F}_q[x]/{\bf F}_q},\phi_f)$ is finite.  Since 
\[ H^n(\Omega^{\updot}_{{\bf F}_q[x]/{\bf F}_q},\phi_f)\simeq {\bf
    F}_q[x_1,\ldots,x_n]/(\partial f/\partial x_1,\ldots,\partial
  f/\partial x_n), \] 
we have
\[ \dim_{{\bf F}_q} H^n(\Omega^{\updot}_{{\bf F}_q[x]/{\bf
    F}_q},\phi_f)=M_f, \]
where $M_f$ denotes the sum of the Milnor numbers of the critical
points of $f$.

We consider one approach to verifying (1.5).  Every $\omega\in\Omega^k_{{\bf
  F}_q[x]/{\bf F}_q}$ can be uniquely written in the form 
\[ \omega=\sum_{1\leq i_1<\cdots<i_k\leq n} \omega(i_1,\ldots,i_k)\,
dx_{i_1}\wedge\cdots\wedge dx_{i_k}, \]
with $\omega(i_1,\ldots,i_k)\in{\bf F}_q[x]$.  If each coefficient
$\omega(i_1,\ldots,i_k)$ is a homogeneous form of degree $l$, we call
$\omega$ {\it homogeneous\/} and write
\begin{align*}
\deg_{\rm coeff}\omega &=l \\
\deg\omega &=l+(n-k)(\delta-1).
\end{align*}
The point of the latter definition is that we can define an
increasing filtration $F.$ on
$\Omega^k_{{\bf F}_q[x]/{\bf F}_q}$ by setting
\[ F_l\Omega^k_{{\bf F}_q[x]/{\bf F}_q} = \text{the ${\bf
  F}_q$-span of homogeneous $k$-forms $\omega$ with $\deg\omega\leq l$}, \]
and $(\Omega^{\updot}_{{\bf F}_q[x]/{\bf
    F}_q},\phi_f)$ then becomes a filtered complex.  Consider the 
associated spectral sequence 
\begin{equation}
E_1^{r,s}=H^{r+s}(F_r/F_{r-1}(\Omega^{\updot}_{{\bf
    F}_q[x]/{\bf F}_q},\phi_f))\Rightarrow
    H^{r+s}(\Omega^{\updot}_{{\bf F}_q[x]/{\bf F}_q},\phi_f).
\end{equation}
The $E_1$-terms are just the cohomology of the homogeneous pieces of
the associated graded complex to $(\Omega^{\updot}_{{\bf
    F}_q[x]/{\bf F}_q},\phi_f)$, which may be identified with
$(\Omega^{\updot}_{{\bf F}_q[x]/{\bf F}_q},\phi_{f^{(\delta)}})$.  Since
this latter complex is isomorphic to the Koszul complex on ${\bf
  F}_q[x_1,\ldots,x_n]$ 
defined by $\{\partial f^{(\delta)}/\partial x_i\}_{i=1}^n$, the assertion
that
\begin{equation}
E_1^{r,s}=0\qquad\text{for $r+s\neq n$}
\end{equation}
is equivalent to the assertion that
\begin{equation}
\{\partial f^{(\delta)}/\partial x_i\}_{i=1}^n\quad\text{form a regular
  sequence in ${\bf F}_q[x_1,\ldots,x_n]$.}
\end{equation}
It follows from (1.8) that
\[ \dim_{{\bf F}_q}{\bf F}_q[x_1,\ldots,x_n]/(\partial
f^{(\delta)}/\partial x_1,\ldots,\partial f^{(\delta)}/\partial
x_n)=(\delta-1)^n, \] 
hence
\begin{equation}
\dim_{{\bf F}_q}\bigoplus_{r+s=n} E_1^{r,s} = (\delta-1)^n. 
\end{equation}
The spectral sequence (1.6) and condition (1.7) imply that (1.5) holds,
and (1.9) then implies that $M_f=(\delta-1)^n$.  Theorem 1.4 of \cite{AS2}, a
slight generalization of Theorem 1.4 above, can then
be reformulated as follows.
\begin{theorem}
Suppose $(1.7)$ holds.  Then $L({\bf A}^n,f;t)^{(-1)^{n+1}}$ is a
polynomial of degree $M_f$, all of whose reciprocal roots have
absolute value $q^{n/2}$.
\end{theorem}

Condition (1.5) holds if there exists a positive integer $e$ such that
\begin{equation}
E_e^{r,s}=0\qquad\text{for $r+s\neq n$}.
\end{equation}
We are interested in determining the extent to which the conclusion of
Theorem~1.10 holds when condition (1.7) is replaced by condition
(1.11) for some $e>1$.  In general, some additional hypothesis
is needed, as is illustrated by the one-variable example
$f(x_1)=x_1^p-x_1$ over the field ${\bf F}_p$.  The purpose of this
paper is to prove a result that provides evidence for such a theorem. 

Dwork has associated to $f$ a complex $(\Omega^{\updot}_{C(b)},D)$ (of
length $n$) depending on a choice of rational parameter $b$ satisfying
$0<b<p/(p-1)$ (we review this theory in section 2).  Each
$\Omega^i_{C(b)}$, $i=0,\ldots,n$, is a $p$-adic Banach space over a
field $\tilde{\Omega}_0$ (a finite extension of ${\bf Q}_p$) and is
equipped with a Frobenius operator $\alpha_i$ commuting with the
differential $D$ of the complex.  Furthermore,
\begin{equation}
L({\bf A}^n,f;t)=\prod_{i=0}^n \det(I-t\alpha_i\mid H^i(
\Omega^{\updot}_{C(b)},D))^{(-1)^{i+1}}. 
\end{equation}

\begin{theorem}
Suppose there exist $e,m$ such that $E_e^{r,s}=0$ for all $r,s$ such
that $r+s=m$.  Then for
\begin{equation}
\frac{\delta}{(p-1)(\delta-e+1)}<b<\frac{p\delta}{(p-1)\delta+e-1}
\end{equation}
we have
\[ H^m(\Omega^{\updot}_{C(b)},D)=0. \]
If $(1.11)$ holds, then, in addition, for $b$ in the range $(1.14)$ we have
\[ \dim_{\tilde{\Omega}_0}H^n(\Omega^{\updot}_{C(b)},D) =M_f. \]
\end{theorem}.

{\it Remark}.  It is easily seen that in (1.14) the upper bound for $b$ is
greater than the lower bound for $b$ if and only if
\begin{equation}
\biggl(1+\frac{p}{(p-1)^2}\biggr)(e-1)<\delta,
\end{equation}
i.~e., (1.15) is equivalent to the existence of rational $b$
satisfying (1.14).  For example, if $e=2$, then this condition
requires $\delta\geq 2$ for odd primes $p$ and $\delta\geq 4$ for
$p=2$.  In general, for $p$ sufficiently large relative to $\delta$,
it becomes simply $e\leq\delta$.

By \cite[section 3.4]{RO}, $\alpha_n$ is invertible on
$H^n(\Omega^{\updot}_{C(b)},D)$, so Theorem 1.13 and equation (1.12)
give the following. 
\begin{corollary}
  Suppose $(1.11)$ holds for a positive integer $e$ satisfying $(1.15)$.
Then $L({\bf A}^n,f;t)^{(-1)^{n+1}}$ is a polynomial of degree $M_f$.
\end{corollary}

The main idea in the proof of Theorem 1.13 is to relate the spectral
sequence (1.6) to the spectral sequence associated to the filtration
by $p$-divisibility on the complex $(\Omega^{\updot}_{C(b)},D)$.
This is accomplished by Theorem 3.1, which is applied in section 4 to
compute the cohomology of $(\Omega^{\updot}_{C(b)},D)$.  

We conjecture that the hypothesis of Corollary 1.16 implies that all
reciprocal roots of $L({\bf A}^n,f;t)^{(-1)^{n+1}}$ have absolute
value $q^{n/2}$.  This would imply, in particular, the estimate
\[ |S({\bf A}^n({\bf F}_{q^i}),f)|\leq M_fq^{ni/2}. \]
Typically, such results are proved by computing the corresponding
$l$-adic cohomology groups.  However, we have been unable to apply our
previous method\cite{AS1,AS2} for calculating $l$-adic cohomology from
$p$-adic cohomology because we have been unable to compute the
$p$-divisibility of the determinant of Frobenius under the hypothesis
of Corollary~1.16.

It is an interesting problem to find geometric conditions that imply
(1.11) for some $e>1$, and we plan to return to this question in a
future article.  As an example, we prove in section 5 the following.
Make (1.3) more precise by writing
\begin{equation}
f=f^{(\delta)}+f^{(\delta')}+f^{(\delta'-1)}+\cdots+f^{(0)},
\end{equation}
where $f^{(j)}$ is homogeneous of degree $j$ and $1\leq \delta'\leq
\delta-1$, i.\ e., $f^{(\delta')}$ is the homogeneous part of
second-highest degree of $f$.  
\begin{theorem}
  Suppose that $f^{(\delta)}=f_1^{a_1}\cdots f_r^{a_r}$, where for every
  subset $\{i_1,\ldots,i_k\} \linebreak \subseteq \{1,\ldots,r\}$ the
  system of equations
\[ f_{i_1}=\cdots=f_{i_k}=0 \]
defines a smooth complete intersection of codimension $k$ in ${\bf
  P}^{n-1}$ and, if either $k\geq 2$ or $k=1$ and $a_{i_1}>1$,  the system
  of equations 
\[ f^{(\delta')}=f_{i_1}=\cdots=f_{i_k}=0 \]
defines a smooth complete intersection of codimension $k+1$ in ${\bf
  P}^{n-1}$.  Suppose also that $(p,\delta\delta'a_1\cdots a_r)=1$.
  Then $(1.11)$ holds for $e=\delta-\delta'+1$.
\end{theorem}

We state here (the proof will appear elsewhere) the first such example
we found.  We refer to Garc\'{\i}a\cite{GA} for the definitions of
``weighted homogeneous'' isolated singularity and ``total degree'' of
a weighted homogeneous isolated singularity.
\begin{theorem}
Suppose that the hypersurface $f^{(\delta)}=0$ in ${\bf P}^{n-1}$ has at
worst weighted homogeneous isolated singularities, of total degrees
$\delta_1,\ldots,\delta_s$, and that none of these singularities lies
on the hypersurface $f^{(\delta')}=0$ in ${\bf P}^{n-1}$.  Suppose also
that $(p,\delta\delta'\delta_1\cdots\delta_s)=1$.  Then $(1.11)$ holds for
$e=\delta-\delta'+1$.
\end{theorem}
The hypothesis of Theorem 1.19 (for $\delta'=\delta-1$) first appears in
Garc\'{\i}a\cite{GA}.  Garc\'{\i}a shows that it implies that the
$l$-adic cohomology groups associated to the exponential sum (1.1)
vanish except in degree $n$, where the cohomology group is pure of
weight $n$ and has dimension $M_f$.  Thus the conclusion of Theorem
1.10 holds in this case.  Garc\'{\i}a's results and Theorem 1.19 are
what originally led us to suspect a result like Theorem 1.13 should
hold.

\section{$p$-adic cohomology}

In this section we review the basic properties of Dwork's $p$-adic
cohomology theory that we shall need.  For a more detailed exposition
of this material, we refer the reader to \cite{AS2}.

Let ${\bf Q}_p$ be the field of $p$-adic numbers,
$\zeta_p$ a primitive $p$-th root of unity, and $\Omega_1={\bf
  Q}_p(\zeta_p)$.  The field $\Omega_1$ is a totally ramified
extension of ${\bf Q}_p$ of degree $p-1$.  Let $K$ be the unramified
extension of ${\bf Q}_p$ of degree $a$.  Set $\Omega_0=K(\zeta_p)$.
The Frobenius automorphism $x\mapsto x^p$ of ${\rm Gal}({\bf F}_q/{\bf
  F}_p)$ lifts to a generator $\tau$ of ${\rm
  Gal}(\Omega_0/\Omega_1)(\simeq {\rm Gal}(K/{\bf Q}_p))$ by requiring
$\tau(\zeta_p)=\zeta_p$.  Let $\Omega$ be the completion of an
algebraic closure of $\Omega_0$.  Denote by ``ord'' the additive
valuation on $\Omega$ normalized by ${\rm ord}\;p=1$ and by ``${\rm
  ord}_q$'' the additive valuation normalized by ${\rm ord}_q\;q=1$.

Let $E(t)$ be the Artin-Hasse exponential series:
\[ E(t)=\exp\biggl(\sum_{i=0}^{\infty} \frac{t^{p^i}}{p^i}\biggr). \]
Let $\gamma\in\Omega_1$ be a solution of $\sum_{i=0}^{\infty}
t^{p^i}/p^i=0$ satisfying ${\rm ord}\;\gamma=1/(p-1)$ and put
\begin{equation}
\theta(t)=E(\gamma t)=\sum_{i=0}^{\infty}
\lambda_it^i\in\Omega_1[[t]].
\end{equation}
The series $\theta(t)$ is a splitting function\cite{DW1} whose
coefficients satisfy 
\begin{equation}
{\rm ord}\;\lambda_i\geq i/(p-1).
\end{equation}

We consider the following spaces of $p$-adic functions.  Let $b$ be a
positive rational number and choose a positive integer $M$ such that
$Mb/p$ and $M\delta/(p(p-1))$ are integers.  Let $\pi$ be such that
\begin{equation}
\pi^{M\delta}=p
\end{equation}
and put $\tilde{\Omega}_1=\Omega_1(\pi)$,
$\tilde{\Omega}_0=\Omega_0(\pi)$.  The element $\pi$ is a uniformizing
parameter for the rings of integers of $\tilde{\Omega}_1$ and
$\tilde{\Omega}_0$.  We extend $\tau\in {\rm Gal}(\Omega_0/\Omega_1)$
to a generator of ${\rm Gal}(\tilde{\Omega}_0/\tilde{\Omega}_1)$ by
requiring $\tau(\pi)=\pi$.  For $u=(u_1,\ldots,u_n)\in{\bf R}^n$, we
put $|u|=u_1+\cdots+u_n$.  Define
\begin{equation}
C(b)=\biggl\{ \sum_{u\in{\bf N}^n} A_u\pi^{Mb|u|}x^u \mid
  \text{$A_u\in\tilde{\Omega}_0$ and $A_u\rightarrow 0$ as
  $u\rightarrow\infty$} \biggr\}. 
\end{equation}
For $\xi=\sum_{u\in{\bf N}^n}A_u\pi^{Mb|u|}x^u\in C(b)$, define
\[ {\rm ord}\;\xi=\min_{u\in{\bf N}^n} \{{\rm ord}\;A_u\}. \]
Given $c\in {\bf R}$, we put
\[ C(b,c)=\{\xi\in C(b)\mid {\rm ord}\;\xi\geq c\}. \]

Let $\hat{f}=\sum_u \hat{a}_ux^u\in K[x_1,\ldots,x_n]$ be the
Teichm\"{u}ller lifting of the polynomial $f\in{\bf
  F}_q[x_1,\ldots,x_n]$, i.~e., $(\hat{a}_u)^q=\hat{a}_u$ and 
the reduction of $\hat{f}$ modulo $p$ is $f$.  Set
\begin{align}
F(x)& = \prod_u \theta(\hat{a}_ux^u), \\
F_0(x)& = \prod_{i=0}^{a-1}\prod_u \theta((\hat{a}_ux^u)^{p^i}).
\end{align}
The estimate (2.2) implies that $F\in
C(b,0)$ for all $b<1/(p-1)$ and $F_0\in C(b,0)$ for all $b<p/q(p-1)$.
Define an operator $\psi$ on formal power series by
\begin{equation}
\psi\biggl(\sum_{u\in{\bf N}^n}A_ux^u\biggr)=\sum_{u\in{\bf N}^n}
A_{pu}x^u.
\end{equation}
It is clear that $\psi(C(b,c))\subseteq C(pb,c)$.
For $0<b<p/(p-1)$, let $\alpha=\psi^a\circ F_0$ be the composition
\[ C(b)\hookrightarrow C(b/q)\xrightarrow{F_0}
C(b/q)\xrightarrow{\psi^a} C(b). \] 
Then $\alpha$ is a completely continuous $\tilde{\Omega}_0$-linear
endomorphism of $C(b)$.  We shall also need to consider
$\beta=\tau^{-1}\circ\psi\circ F$, which is a completely continuous
$\tilde{\Omega}_1$-linear (or $\tilde{\Omega}_0$-semilinear)
endomorphism of $C(b)$.  Note that $\alpha=\beta^a$.

Set $\hat{f}_i=\partial\hat{f}/\partial x_i$ and let
$\gamma_l=\sum_{i=0}^l \gamma^{p^i}/p^i$.  By the definition of
$\gamma$, we have 
\begin{equation}
{\rm ord}\;\gamma_l\geq \frac{p^{l+1}}{p-1}-l-1.
\end{equation}
For $i=1,\ldots,n$, define differential operators $D_i$ by
\begin{equation}
D_i=\pi^{\epsilon}\frac{\partial}{\partial x_i}+\pi^{\epsilon}H_i,
\end{equation}
where $\epsilon=Mb(\delta-1)-M\delta/(p-1)$ and 
\begin{equation}
H_i=\sum_{l=0}^{\infty} \gamma_lp^lx_i^{p^l-1}\hat{f}_i^{\tau^l}(x^{p^l}) \in
C\biggl(b,\frac{1}{p-1}-b\frac{\delta-1}{\delta}\biggr)
\end{equation}
for $b<p/(p-1)$.  Thus $D_i$ and ``multiplication by $H_i$'' operate
on $C(b)$ for $b<p/(p-1)$.  As explained in \cite{AS2}, we have
\begin{align}
\alpha\circ x_iD_i&=qx_iD_i\circ\alpha, \\
\beta\circ x_iD_i&=px_iD_i\circ\beta.
\end{align}
The significance of the normalizing factor $\pi^{\epsilon}$ will be
explained below.

Consider the de Rham-type complex $(\Omega_{C(b)}^{\updot},D)$, where
\[ \Omega_{C(b)}^k=\bigoplus_{1\leq i_1<\cdots<i_k\leq n} C(b)\,
dx_{i_1}\wedge\cdots\wedge dx_{i_k} \]
and $D:\Omega_{C(b)}^k\rightarrow \Omega_{C(b)}^{k+1}$ is defined by
\[ D(\xi\,dx_{i_1}\wedge\cdots\wedge dx_{i_k})=\biggl(\sum_{i=1}^n
D_i(\xi)\,dx_i\biggr)\wedge dx_{i_1}\wedge\cdots\wedge dx_{i_k}. \]
We extend the mapping $\alpha$ to a mapping
$\alpha_{\updot}:\Omega_{C(b)}^{\updot}\rightarrow
\Omega_{C(b)}^{\updot}$ defined by linearity and the formula
\[ \alpha_k(\xi\,dx_{i_1}\wedge\cdots\wedge dx_{i_k})=
q^{n-k}\frac{1}{x_{i_1}\cdots x_{i_k}}\alpha(x_{i_1}\cdots
x_{i_k}\xi)\, dx_{i_1}\wedge\cdots\wedge dx_{i_k}. \]
Equation (2.11) implies that $\alpha_{\updot}$ is a map of complexes.
The Dwork trace formula, as formulated by Robba\cite{RO}, then gives
\begin{equation}
L({\bf A}^n,f;t)=\prod_{k=0}^n \det(I-t\alpha_k\mid
\Omega_{C(b)}^k)^{(-1)^{k+1}}.
\end{equation}
This implies
\begin{equation}
L({\bf A}^n,f;t)=\prod_{k=0}^n \det(I-t\alpha_k\mid
H^k(\Omega_{C(b)}^{\updot},D))^{(-1)^{k+1}},
\end{equation}
where we denote the induced map on cohomology by $\alpha_k$ also.

The $p$-adic Banach space $C(b)$ has a decreasing filtration
$\{\hat{F}^rC(b)\}_{r=-\infty}^{\infty}$ defined by setting
\[ \hat{F}^rC(b)=\bigg\{\sum_{u\in{\bf N}^n} A_u\pi^{Mb|u|}x^u\in
C(b) \mid A_u\in\pi^r{\mathcal O}_{\tilde{\Omega}_0} \text{ for all
  $u$}\bigg\}, \]
where ${\mathcal O}_{\tilde{\Omega}_0}$ denotes the ring of integers
of $\tilde{\Omega}_0$.  We extend this to a filtration on
$\Omega^{\updot}_{C(b)}$ by defining
\[ \hat{F}^r\Omega^k_{C(b)}=\bigoplus_{1\leq i_1<\cdots<i_k\leq n}
  \hat{F}^rC(b)\, dx_{i_1}\wedge\cdots\wedge dx_{i_k}. \]
This filtration is exhaustive and separated, i.~e.,
\[ \bigcup_{r\in{\bf Z}}
\hat{F}^r\Omega^{\updot}_{C(b)}=\Omega^{\updot}_{C(b)} \quad\text{and}\quad
\bigcap_{r\in{\bf Z}}\hat{F}^r\Omega^{\updot}_{C(b)}=(0). \]
Our choice of the normalizing factor $\pi^{\epsilon}$ in (2.9)
guarantees that the $D_i$ respect this filtration, i.~e.,
$D_i(\hat{F}^rC(b))\subseteq \hat{F}^rC(b)$, hence
$D(\hat{F}^r\Omega^k_{C(b)})\subseteq\hat{F}^r\Omega^{k+1}_{C(b)}$.
Associated to the filtered complex
$(\Omega^{\updot}_{C(b)},D)$ is the spectral sequence
\begin{equation}
\hat{E}_1^{r,s}=H^{r+s}(\hat{F}^r\Omega^{\updot}_{C(b)}/
\hat{F}^{r+1}\Omega^{\updot}_{C(b)})\Rightarrow
H^{r+s}(\Omega^{\updot}_{C(b)},D).  
\end{equation}
The notation $\hat{E}^{r,s}_t$ does not express the dependence of this
spectral sequence on the choice of $b$, however, this should not cause
confusion.  We shall prove Theorem 1.13 by analyzing this spectral
sequence. 

For notational convenience we define an ``exterior derivative''
$d:\Omega^k_{C(b)}\rightarrow \Omega^{k+1}_{C(b)}$.  It is
characterised by $\tilde{\Omega}_0$-linearity and the formula
\[ d(\xi\,dx_{i_1}\wedge\cdots\wedge dx_{i_k})=\biggl(\sum_{i=1}^n
\frac{\partial \xi}{\partial x_i}\,dx_i\biggr)\wedge
dx_{i_1}\wedge\cdots\wedge dx_{i_k}. \] 
Although we use the same symbol ``$d$'' for the exterior derivative on
both $\Omega^{\updot}_{C(b)}$ and $\Omega^{\updot}_{{\bf F}_q[x]/{\bf
    F}_q}$, its meaning will be clear from the context.  

It will also be convenient to define a ``reduction map'' ${\rm
  red}:\hat{F}^0\Omega^m_{C(b)}\rightarrow\Omega^m_{{\bf F}_q[x]/{\bf
    F}_q}$ as follows.  For
\[ \omega=\sum_{1\leq i_1<\cdots<i_m\leq n}\biggl(\sum_{u\in{\bf N}^n}
  A_u(i_1,\ldots,i_m)\pi^{Mb|u|}x^u\biggr)\,dx_{i_1}\wedge\cdots
  \wedge dx_{i_m}\in\hat{F}^0\Omega^m_{C(b)}, \]
define
\begin{equation}
{\rm red}(\omega)=\sum_{1\leq i_1<\cdots<i_m\leq n}\biggl(\sum_{u\in{\bf N}^n}
  \bar{A}_u(i_1,\ldots,i_m)x^u\biggr)\,dx_{i_1}\wedge\cdots
  \wedge dx_{i_m}\in\Omega^m_{{\bf F}_q[x]/{\bf F}_q},
\end{equation}
where $\bar{A}_u(i_1,\ldots,i_m)\in{\bf F}_q$ is the reduction mod
$\pi$ of $A_u(i_1,\ldots,i_m)\in{\mathcal O}_{\tilde{\Omega}_0}$.
Note that the inner sum in (2.16) is finite by the definition of
$C(b)$ (equation (2.4)).

\section{Relation between the spectral sequences $E^{r,s}_t$ and
$\hat{E}^{r,s}_t$} 

The key to our computation of $p$-adic cohomology is the
following.
\begin{theorem}
Fix a positive integer $e$, let
$m\in\{0,1,\ldots,n\}$, and let $b$ be a rational number satisfying
$(1.14)$.  Let $\{\phi_i\}_{i\in I}\subseteq\Omega^m_{{\bf
    F}_q[x]/{\bf F}_q}$, $I$ an index set, be a set of homogeneous
$m$-forms such that for each $r\geq 0$, the classes
$\{[\phi_i]\}_{\deg\phi_i=r}$ form a basis for $E_e^{r,m-r}$ as ${\bf
  F}_q$-vector space.  Let $\{\tilde{\phi}_i\}_{i\in I}\subseteq
\hat{F}^0\Omega^m_{C(b)}$ be a set of $m$-forms satisfying ${\rm
  red}(\tilde{\phi}_i)=\phi_i$ for all $i\in I$.  Then for each
$r\in{\bf Z}$, the classes $\{[\pi^r\tilde{\phi}_i]\}_{i\in I}$ form a
basis for $\hat{E}_{Mb(e-1)+1}^{r,m-r}$ as ${\bf F}_q$-vector space.
\end{theorem}

{\it Remark}.  We observe for later reference that (1.14) implies ${\rm
  ord}_p\;\pi^{Mb(e-1)}<1$.  

{\it Proof}.  Consider the operator
$D_i=\pi^{\epsilon}\partial/\partial x_i+\pi^{\epsilon}H_i$.  The
terms of degree $<\delta-e$ in $\pi^{\epsilon}H_i$ lie in
$\hat{F}^{Mbe}C(b)$.  The upper bound on $b$ given by (1.14) guarantees
that the terms of degree $>\delta-1$ in $\pi^{\epsilon}H_i$ lie in
$\hat{F}^{Mb(e-1)+1}C(b)$.  The lower bound on $b$ given by (1.14)
guarantees that $\pi^{\epsilon}\partial/\partial x_i$ maps
$\hat{F}^0C(b)$ into $\hat{F}^{Mb(e-1)+1}C(b)$.  We thus have for
$\xi\in \hat{F}^0C(b)$,
\begin{equation}
D_i(\xi)-\pi^{Mb(\delta-1)}\biggl(\sum_{j=0}^{e-1}\frac{\partial
  \hat{f}^{(\delta-j)}}{\partial x_i}\biggr)\xi\in \hat{F}^{Mb(e-1)+1}C(b).
\end{equation}
It follows that if $\omega\in \hat{F}^r\Omega^k_{C(b)}$, $0\leq k\leq
n$, then
\begin{equation}
D(\omega)-\pi^{Mb(\delta-1)}\biggl( \sum_{j=0}^{e-1}
d\hat{f}^{(\delta-j)}\biggr)\wedge\omega\in
\hat{F}^{r+Mb(e-1)+1}\Omega^k_{C(b)}. 
\end{equation}
Note that
\begin{equation}
\pi^{Mb(\delta-j-1)}d\hat{f}^{(\delta-j)}\in \hat{F}^0\Omega^1_{C(b)}
\qquad\text{for $j=0,1,\ldots,e-1$}. 
\end{equation}

We recall the definition of $\hat{E}^{r,s}_t$.  Put
\[ Z_t^{r,s}=\{\omega\in \hat{F}^r\Omega^{r+s}_{C(b)}\mid
D(\omega)\in \hat{F}^{r+t}\Omega^{r+s+1}_{C(b)}\}. \]
Then
\[ \hat{E}^{r,s}_t=\frac{Z_t^{r,s}+\hat{F}^{r+1}\Omega^{r+s}_{C(b)}}
{D(Z_{t-1}^{r-t+1,s+t-2})+\hat{F}^{r+1}\Omega^{r+s}_{C(b)}}.
\]

Theorem 3.1 asserts that if $\omega\in
Z_{Mb(e-1)+1}^{r,m-r}$, then there exist a finite subset $I_0\subseteq
I$, a collection $\{c_i\}_{i\in I_0}\subseteq{\mathcal
  O}_{\tilde{\Omega}_0}$ uniquely determined mod $\pi$, and
\[ \xi\in Z_{Mb(e-1)}^{r-Mb(e-1),m-r-1+Mb(e-1)} \]
such that
\[ \omega\equiv \sum_{i\in I_0} c_i\pi^r\tilde{\phi}_i+ D(\xi)
\pmod{\hat{F}^{r+1}\Omega^{m}_{C(b)}}. \]
Since $\pi^r\hat{F}^0=\hat{F}^r$, we may reduce to the case $r=0$ by
multiplication by $\pi^{-r}$.  So let
$\omega\in\hat{F}^0\Omega^m_{C(b)}$ be such that 
\begin{equation}
D(\omega)\in \hat{F}^{Mb(e-1)+1}\Omega^{m+1}_{C(b)}.
\end{equation}
We must show that there exist $\{c_i\}_{i\in I_0}\subseteq{\mathcal
  O}_{\tilde{\Omega}_0}$ uniquely determined mod $\pi$ and $\xi\in
  \hat{F}^{-Mb(e-1)}\Omega^{m-1}_{C(b)}$ such that 
\begin{equation}
\omega\equiv \sum_{i\in I_0} c_i\tilde{\phi}_i + D(\xi)
\pmod{\hat{F}^1\Omega^m_{C(b)}}.
\end{equation}
In view of (3.3), (3.5) is equivalent to the condition
\begin{equation}
\biggl(\sum_{j=0}^{e-1}\pi^{Mb(\delta-1)}d\hat{f}^{(\delta-j)}\biggr)
\wedge\omega\in \hat{F}^{Mb(e-1)+1}\Omega^{m+1}_{C(b)} 
\end{equation}
and (3.6) is equivalent to the condition
\begin{equation}
\omega\equiv \sum_{i\in I_0} c_i\tilde{\phi}_i + 
\biggl(\sum_{j=0}^{e-1}\pi^{Mb(\delta-1)}d\hat{f}^{(\delta-j)}\biggr)\wedge\xi
\pmod{\hat{F}^1\Omega^m_{C(b)}}. 
\end{equation}
This reduces the proof of Theorem 3.1 to showing that, given
$\omega\in \hat{F}^0\Omega^m_{C(b)}$ satisfying
(3.7), there exist a finite set $\{c_i\}_{i\in I_0}\subseteq{\mathcal
  O}_{\tilde{\Omega}_0}$, uniquely determined mod $\pi$, and $\xi\in
\hat{F}^{-Mb(e-1)}\Omega_{C(b)}^{m-1}$ satisfying (3.8).  

Equations (3.4) and (3.7) imply that for $i=0,1,\ldots,e-2$ we have
\[ \biggl(\sum_{j=0}^i \pi^{Mb(\delta-1)}d\hat{f}^{(\delta-j)}\biggr)
\wedge\omega\in \hat{F}^{Mb(i+1)}\Omega^{m+1}_{C(b)} . \]
We may thus define $\eta_i\in \hat{F}^0\Omega^{m+1}_{C(b)}$ for
$i=-1,0,1,\ldots,e-2$ recursively as follows.  Set $\eta_{-1}=0$ and
for $i=0,1,\ldots,e-2$ define $\eta_i$ by the formula
\begin{equation}
\eta_{i-1}+\pi^{Mb(\delta-1-i)}d\hat{f}^{(\delta-i)}
\wedge\omega=\pi^{Mb}\eta_i. 
\end{equation}
From (3.7) we also get that
\begin{equation}
\eta_{e-2}+\pi^{Mb(\delta-e)}d\hat{f}^{(\delta-e+1)}\wedge\omega \equiv 0
\pmod{\hat{F}^1\Omega^{m+1}_{C(b)}}.
\end{equation}

Write
\begin{align}
\omega &=\sum_{k=0}^{\infty} \pi^{Mbk}\omega^{(k)} \\
\eta_i &=\sum_{k=0}^{\infty} \pi^{Mbk}\eta_i^{(k)}, 
\end{align}
where $\omega^{(k)}$ (resp.\ $\eta_i^{(k)}$) is an $m$-form (resp.\
$(m+1)$-form) whose coefficients are homogeneous polynomials of degree
$k$ with coefficients in ${\mathcal O}_{\tilde{\Omega}_0}$.
Substituting (3.11) and (3.12) into (3.9) and cancelling a power of
$\pi$ gives 
\begin{equation}
\eta_{i-1}^{(k+\delta-1-i)} + d\hat{f}^{(\delta-i)}\wedge\omega^{(k)} =
\pi^{Mb}\eta_i^{(k+\delta-1-i)} 
\end{equation}
for $i=0,1,\ldots,e-2$ and all $k$.  From (3.10) we also have
\begin{equation}
\eta_{e-2}^{(k+\delta-1-i)}+d\hat{f}^{(\delta-e+1)}\wedge \omega^{(k)}
\equiv 0 \pmod{\pi}.
\end{equation}

Since $\tilde{\Omega}_0=K(\pi)$, where $K$ is an unramified extension
of ${\bf Q}_p$ and $\pi$ is an $(M\delta)$-th root of $p$, we may write
\begin{equation}
\omega^{(k)}=\sum_{l=0}^{\infty} \pi^l\omega^{(k)}_l,
\end{equation}
where $\omega^{(k)}_l$ is an $m$-form whose coefficients are
homogeneous polynomials with coefficients in ${\mathcal O}_K$, the
ring of integers of $K$.  For $i=0,1,\ldots,e-2$ and $k\geq 0$, put
\begin{equation}
\zeta^{(k)}_i=\sum_{l=0}^{\infty} \pi^l\omega^{(k)}_{l+Mb(i+1)}.
\end{equation}

\begin{lemma}
For $i=0,1,\ldots,e-2$ and $k\geq 0$, we have
\begin{equation}
\eta_i^{(k+\delta-1-i)}\equiv \sum_{j=0}^i
d\hat{f}^{(\delta-j)}\wedge\zeta_{i-j}^{(k-i+j)}
\pmod{\pi^{Mb(e-i-2)+1}}.
\end{equation}
\end{lemma}

{\it Remark}.  The congruence (3.18) is to be interpreted as meaning
that the difference between the two sides is an $(m+1)$-form all of
whose coefficients are homogeneous polynomials of degree $k+\delta-1-i$
with coefficients in $\pi^{Mb(e-i-2)+1}{\mathcal
  O}_{\tilde{\Omega}_0}$.

{\it Proof}.  The assertion is vacuous for $e=1$, so we assume $e\geq
2$.  We fix $k+\delta-1-i$ and prove (3.18) by induction on $i$.
By (3.13) with $i=0$ and (3.15) we have
\begin{equation}
d\hat{f}^{(\delta)}\wedge\sum_{l=0}^{\infty} \pi^l\omega^{(k)}_l =
\pi^{Mb}\eta_0^{(k+\delta-1)},
\end{equation}
which implies
\[ d\hat{f}^{(\delta)}\wedge\omega^{(k)}_0\equiv 0 \pmod{\pi}. \]
But since the left-hand side has coefficients in ${\mathcal O}_K$,
this congruence must actually hold mod $p$.  Suppose we have proved
\begin{equation}
d\hat{f}^{(\delta)}\wedge\omega^{(k)}_l\equiv 0 \pmod{p}
\end{equation}
for $l=0,1,\ldots,L$, where $0\leq L<Mb-1$.  Since
$\pi^{Mb(e-1)+1}$ divides $p$ (see the remark following Theorem 3.1)
it follows from (3.19) that
\[ d\hat{f}^{(\delta)}\wedge\sum_{l=L+1}^{\infty}
\pi^l\omega_l^{(k)} \equiv 0\pmod{\pi^{Mb}}. \]
Dividing by $\pi^{L+1}$ gives
\[ d\hat{f}^{(\delta)}\wedge\omega^{(k)}_{L+1}\equiv 0\pmod{\pi}, \]
and again the left-hand side has coefficients in ${\mathcal O}_K$, so
this congruence holds mod~$p$.  This proves that (3.20) holds for $0\leq
l\leq Mb-1$.  Equation (3.19) then implies
\[ d\hat{f}^{(\delta)}\wedge\sum_{l=Mb}^{\infty}\pi^l\omega^{(k)}_l
\equiv \pi^{Mb}\eta_0^{(k+\delta-1)} \pmod{\pi^{Mb(e-1)+1}}. \]
Dividing this congruence by $\pi^{Mb}$ gives (3.18) for $i=0$.

Suppose (3.18) holds for some $i$, $0\leq i<e-2$.  We prove (3.18) for
$i+1$.  From (3.13) and the induction hypothesis we have
\begin{multline*}
\sum_{j=0}^i d\hat{f}^{(\delta-j)}\wedge\zeta_{i-j}^{(k+j-i)} +
d\hat{f}^{(\delta-i-1)}\wedge\omega^{(k+1)}\equiv \\
\pi^{Mb}\eta_{i+1}^{(k+\delta-1-i)} \pmod{\pi^{Mb(e-i-2)+1}}. 
\end{multline*}
Using (3.15) and (3.16), this is equivalent to
\begin{multline}
\sum_{j=0}^i d\hat{f}^{(\delta-j)}\wedge\biggl(\sum_{l=0}^{\infty}
\pi^l\omega_{l+Mb(i-j+1)}^{(k+j-i)}\biggr) +
d\hat{f}^{(\delta-i-1)}\wedge \biggl(\sum_{l=0}^{\infty}
\pi^l\omega_l^{(k+1)}\biggr) \equiv \\
\pi^{Mb}\eta_{i+1}^{(k+\delta-1-i)}\pmod{\pi^{Mb(e-i-2)+1}}. 
\end{multline}
Since $e-i>2$, this implies
\[ \sum_{j=0}^i
d\hat{f}^{(\delta-j)}\wedge\omega^{(k+j-i)}_{Mb(i-j+1)} +
d\hat{f}^{(\delta-i-1)}\wedge\omega^{(k+1)}_0\equiv 0\pmod{\pi}, \]
and since the left-hand side has coefficients in ${\mathcal O}_K$,
this congruence holds mod $p$.  Arguing by induction on $l$ exactly as
in the case $i=0$ gives
\[ \sum_{j=0}^i
d\hat{f}^{(\delta-j)}\wedge\omega^{(k+j-i)}_{l+Mb(i-j+1)} +
d\hat{f}^{(\delta-i-1)}\wedge\omega^{(k+1)}_l\equiv 0\pmod{p} \]
for $0\leq l\leq Mb-1$.  Equation (3.21) then implies
\begin{multline*}
\sum_{j=0}^i d\hat{f}^{(\delta-j)}\wedge\biggl(\sum_{l=Mb}^{\infty}
\pi^l\omega^{(k+j-i)}_{l+Mb(i-j+1)}\biggr) +
d\hat{f}^{(\delta-i-1)}\wedge\biggl(\sum_{l=Mb}^{\infty}
\pi^l\omega^{(k+1)}_l\biggr)\equiv \\
\pi^{Mb}\eta_{i+1}^{(k+\delta-1-i)}\pmod{\pi^{Mb(e-i-2)+1}}. 
\end{multline*}
Dividing by $\pi^{Mb}$ now gives
\begin{multline*}
\eta_{i+1}^{(k+\delta-1-i)} \equiv \sum_{j=0}^i
d\hat{f}^{(\delta-j)}\wedge \biggl(\sum_{l=0}^{\infty}
\pi^l\omega_{l+Mb(i-j+2)}^{(k+j-i)}\biggr) +
d\hat{f}^{(\delta-i-1)}\wedge\biggl(\sum_{l=0}^{\infty}
\pi^l\omega^{(k+1)}_{l+Mb}\biggr) \\
\pmod{\pi^{Mb(e-i-3)+1}}. 
\end{multline*}
Taking into account (3.16), we see that this is just (3.18) for $i+1$.

Combining Lemma 3.17 with equations (3.13) and (3.14) gives
\begin{equation}
\sum_{j=0}^{i-1} d\hat{f}^{(\delta-j)}\wedge\zeta_{i-1-j}^{(k+j-i)}
+ d\hat{f}^{(\delta-i)}\wedge\omega^{(k)}\equiv 0\pmod{\pi}
\end{equation}
for $i=0,1,\ldots,e-1$.

For $i\in I$, write
\[ \tilde{\phi}_i=\sum_{k=0}^{\infty} \pi^{Mbk}\tilde{\phi}_i^{(k)},
\]
where $\tilde{\phi}_i^{(k)}$ is an $m$-form whose coefficients are
homogeneous polynomials of degree $k$ with coefficients in ${\mathcal
  O}_{\tilde{\Omega}_0}$.  Our hypothesis that ${\rm
  red}(\tilde{\phi}_i)=\phi_i$ implies that if $\deg_{\rm
  coeff}\phi_i=k_0$, then $\phi_i$ is the reduction mod $\pi$ of
$\tilde{\phi}_i^{(k_0)}$, while for $k'\neq k_0$, the reduction mod
$\pi$ of $\tilde{\phi}_i^{(k')}$ is 0.  

\begin{lemma}
For each $k$ there exist $\{c_i\}_{\deg_{\rm coeff}\phi_i=k}\subseteq
{\mathcal O}_K$ and
$\{\xi_j^{(k-\delta+j)}\}_{j=1}^e\subseteq\Omega^{m-1}_{C(b)}$, where
$\xi_j^{(k-\delta+j)}$ is an $(m-1)$-form whose coefficients are
homogeneous polynomials of degree $k-\delta+j$ with 
coefficients in ${\mathcal O}_K$, such that
\begin{equation}
\omega^{(k)}\equiv \sum_{\deg_{\rm coeff}\phi_i=k}
c_i\tilde{\phi}_i^{(k)} + \sum_{j=0}^{e-1} d\hat{f}^{(\delta-j)}\wedge
\xi_{j+1}^{(k-\delta+j+1)} \pmod{\pi}
\end{equation}
and such that
\begin{equation}
\sum_{j=0}^i d\hat{f}^{(\delta-j)}\wedge\xi_{j+e-i}^{(k-\delta+j+e-i)}
\equiv 0 \pmod{p}
\end{equation}
for $i=0,1,\ldots,e-2$.
\end{lemma}

{\it Remark}.  By the definition of $C(b)$, there exists a positive
integer $k_{\omega}$ such that $\omega^{(k)}\equiv 0\pmod{\pi}$ for
all $k>k_{\omega}$.  For $k>k_{\omega}$ we take all $c_i$ and all
$\xi_j^{(k-\delta+j)}$ equal to 0, which satisfies the conclusion of the
lemma.

{\it Proof}.  Let $\bar{\omega}^{(k)}\in\Omega^m_{{\bf
    F}_q[x]/{\bf F}_q}$ (resp.\ $\bar{\zeta}_j^{(k)}\in
\Omega^m_{{\bf F}_q[x]/{\bf F}_q}$) be the reduction mod $\pi$ of
$\omega^{(k)}$ (resp.\ $\zeta_j^{(k)}$).  Then (3.22) implies
\[ \sum_{j=0}^{i-1} df^{(\delta-j)}\wedge\bar{\zeta}_{i-1-j}^{(k-i+j)} +
    df^{(\delta-i)}\wedge\bar{\omega}^{(k)}=0 \]
for $i=0,1,\ldots,e-1$.  The definition of $E^{r,m-r}_e$ and the
hypothesis that $\{\phi_i\}_{\deg\phi_i=r}$ spans
$E^{r,m-r}_e$ imply that there exist $\{\bar{c}_i\}_{\deg_{\rm
    coeff}\phi_i=k}\subseteq {\bf F}_q$ and
$\{\bar{\xi}_j^{(k-\delta+j)}\}_{j=1}^e \subseteq \Omega^{m-1}_{{\bf
    F}_q[x]/{\bf F}_q}$ such that 
\begin{equation}
\bar{\omega}^{(k)}=\sum_{\deg_{\rm coeff}\phi_i=k} \bar{c}_i\phi_i +
\sum_{j=0}^{e-1} df^{(\delta-j)}\wedge\bar{\xi}_{j+1}^{(k-\delta+j+1)}  
\end{equation}
and such that
\begin{equation}
\sum_{j=0}^i df^{(\delta-j)}\wedge\bar{\xi}_{j+e-i}^{(k+j+e-i-\delta)} = 0
\end{equation}
for $i=0,1,\ldots,e-2$.  

Let $c_i\in{\mathcal O}_K$ be any lifting of $\bar{c}_i\in{\bf F}_q$
and let $\xi_j^{(k+j-\delta)}\in\Omega^{m-1}_{C(b)}$ be any lifting of 
$\bar{\xi}_j^{(k+j-\delta)}$ that has coefficients in ${\mathcal
  O}_K$.  Thus  $\xi_j^{(k+j-\delta)}$ is an $(m-1)$-form whose
coefficients are homogeneous polynomials of degree $k+j-\delta$ with
coefficients in~${\mathcal O}_K$.  Then (3.26) and (3.27) imply (3.24)
and (3.25), the congruence (3.25) holding mod~$p$ rather than just mod
$\pi$ because the left-hand side has coefficients in ${\mathcal
  O}_K$.  This completes the proof of Lemma 3.23.

We continue with the proof of Theorem 3.1.  Put
\[ \xi_i=\sum_{k=0}^{\infty} \pi^{Mb(k-\delta+i)}\xi_i^{(k-\delta+i)}\in
\hat{F}^0\Omega^{m-1}_{C(b)} \]
for $i=1,2,\ldots,e$.  Let $I_0\subseteq I$ be the (finite) subset
consisting of all indices $i$ such that $\deg_{\rm coeff}\phi_i\leq
k_{\omega}$.  Then (3.11), (3.24), and the remark following Lemma 3.23
imply 
\begin{equation}
\omega\equiv \sum_{i\in I_0}c_i\tilde{\phi}_i + \sum_{j=0}^{e-1}
\pi^{Mb(\delta-j-1)}d\hat{f}^{(\delta-j)}\wedge\xi_{j+1}
\pmod{\hat{F}^1\Omega^m_{C(b)}}. 
\end{equation}
Since $p$ is divisible by $\pi^{Mb(e-1)+1}$, (3.25) implies 
\begin{equation}
\sum_{j=0}^i \pi^{Mb(\delta-j-1)}d\hat{f}^{(\delta-j)}
\wedge\xi_{j+e-i}\equiv 0 \pmod{\hat{F}^{Mb(e-1)+1}\Omega^m_{C(b)}}
\end{equation}
for $i=0,1,\ldots,e-2$.  Put
\begin{equation}
\xi=\sum_{j=0}^{e-1} \pi^{-Mbj}\xi_{j+1}\in
\hat{F}^{-Mb(e-1)}\Omega^{m-1}_{C(b)}. 
\end{equation}
With this choice of $\xi$, the right-hand side of (3.8) becomes
\begin{equation}
\sum_{i\in I_0} c_i\tilde{\phi}_i + \biggl(\sum_{i=0}^{e-1}
\pi^{Mbi}(\pi^{Mb(\delta-i-1)}d\hat{f}^{(\delta-i)})\biggr)\wedge
\biggl(\sum_{j=0}^{e-1} \pi^{-Mbj}\xi_{j+1}\biggr). 
\end{equation}
When (3.31) is expanded, the wedge product of a pair of terms with $i>j$
lies in $\hat{F}^{Mb}\Omega^m_{C(b)}$.  Putting $k=j-i$ when $j\geq i$, we
see that (3.31) is congruent mod $\hat{F}^{Mb}\Omega^m_{C(b)}$ to
\begin{equation}
\sum_{i\in I_0}c_i\tilde{\phi}_i+\sum_{k=0}^{e-1}\pi^{-Mbk} \sum_{l=0}^{e-1-k}
\pi^{Mb(\delta-l-1)}d\hat{f}^{(\delta-l)}\wedge\xi_{l+k+1}. 
\end{equation}
For $k=1,\ldots,e-1$, the corresponding summand of (3.32) lies in
$\hat{F}^1\Omega^m_{C(b)}$ by (3.29), hence (3.32) is congruent mod
$\hat{F}^1\Omega^m_{C(b)}$ to
\[ \sum_{i\in I_0}c_i\tilde{\phi}_i + \sum_{l=0}^{e-1} \pi^{Mb(\delta-l-1)}
d\hat{f}^{(\delta-l)}\wedge\xi_{l+1}. \]
But this is $\equiv\omega\pmod{\hat{F}^1\Omega^m_{C(b)}}$ by (3.28).  
Thus congruence (3.8) holds.  

To complete the proof of Theorem 3.1, it remains to show that the
$c_i$ are uniquely determined mod $\pi$.  Suppose there exist a finite
subset $I_0\subseteq I$, $\{c_i\}_{i\in I_0}\subseteq {\mathcal
  O}_{\tilde{\Omega}_0}$ and
$\xi\in\hat{F}^{-Mb(e-1)}\Omega^{m-1}_{C(b)}$ such that 
\begin{equation}
\sum_{i\in I_0} c_i\tilde{\phi}_i \equiv
\sum_{j=0}^{e-1}\pi^{Mb(\delta-1)}d\hat{f}^{(\delta-j)}\wedge\xi
\pmod{\hat{F}^1\Omega^m_{C(b)}}.
\end{equation}
We must show $c_i\equiv 0\pmod{\pi}$ for all $i\in I_0$.  Write
\begin{equation}
\xi=\pi^{-Mb(e-1)}\sum_{k=0}^{\infty} \pi^{Mbk}\xi^{(k)}, 
\end{equation}
where the coefficients of $\xi^{(k)}\in\Omega^{m-1}_{C(b)}$ are
homogeneous forms of degree $k$ with coefficients in ${\mathcal
  O}_{\tilde{\Omega}_0}$.  Substituting (3.34) into (3.33) and
cancelling $\pi^{Mbk}$ gives
\begin{equation}
\sum_{\substack{i\in I_0\\ \deg_{\rm coeff}\phi_i=k}}
c_i\tilde{\phi}_i^{(k)} \equiv 
\sum_{j=0}^{e-1}\pi^{Mb(j-e+1)}d\hat{f}^{(\delta-j)}\wedge\xi^{(k-\delta+j+1)}
\pmod{\pi}
\end{equation}
for all $k$.  We may thus define $\eta_i\in\Omega^{m-1}_{C(b)}$ for
$i=-1,0,1,\ldots,e-2$ recursively as follows.  Put $\eta_{-1}=0$ and
define $\eta_i$ for $i=0,1,\ldots,e-2$ by 
\begin{equation}
\eta^{(k)}_{i-1}+d\hat{f}^{(\delta-i)}\wedge\xi^{(k-\delta+i+1)}
=\pi^{Mb}\eta_i^{(k)}. 
\end{equation}
We also have from (3.35) that
\begin{equation}
\eta^{(k)}_{e-2}+d\hat{f}^{(\delta-e+1)}\wedge\xi^{(k-\delta+e)} \equiv
\sum_{\substack{i\in I_0\\ \deg_{\rm coeff}\phi_i=k}}
c_i\tilde{\phi}_i^{(k)} \pmod{\pi}. 
\end{equation}
By a shift of index, (3.36) and (3.37) become
\begin{equation}
\eta_{i-1}^{(k+\delta-1-i)}+d\hat{f}^{(\delta-i)}\wedge\xi^{(k)} =
\pi^{Mb}\eta_i^{(k+\delta-1-i)} 
\end{equation}
for $i=0,1,\ldots,e-2$ and
\begin{equation}
\eta_{e-2}^{(k+\delta-e)}+d\hat{f}^{(\delta-e+1)}\wedge\xi^{(k)} \equiv
\sum_{\substack{i\in I_0\\ \deg_{\rm coeff}\phi_i=k+\delta-e}}
c_i\tilde{\phi}_i^{(k+\delta-e)} 
\pmod{\pi}. 
\end{equation}

Write $\xi^{(k)}=\sum_{l=0}^{\infty} \pi^l\xi_l^{(k)}$, where
$\xi_l^{(k)}$ is an $(m-1)$-form whose coefficients are homogeneous
polynomials of degree $k$ with coefficients in ${\mathcal O}_K$.  Put
\[ \rho_i^{(k)}=\sum_{l=0}^{\infty} \pi^l\xi^{(k)}_{l+Mb(i+1)}. \]
Since (3.38) is identical in form to (3.13), we may apply Lemma 3.17
to conclude that
\begin{equation}
\eta_i^{(k+\delta-1-i)}\equiv \sum_{j=0}^i
d\hat{f}^{(\delta-j)}\wedge\rho_{i-j}^{(k-i+j)} \pmod{\pi^{Mb(e-i-2)+1}}.
\end{equation}
Substitution into (3.38) and (3.39) then gives
\begin{equation}
\sum_{j=0}^{i-1}d\hat{f}^{(\delta-j)}\wedge\rho_{i-1-j}^{(k+j-i)} +
d\hat{f}^{(\delta-i)}\wedge\xi^{(k)}\equiv 0\pmod{\pi}
\end{equation}
for $i=0,1,\ldots,e-2$ and
\begin{multline}
\sum_{j=0}^{e-2} d\hat{f}^{(\delta-j)}\wedge\rho_{e-2-j}^{(k+j-e-1)} +
d\hat{f}^{(\delta-e+1)}\wedge\xi^{(k)} \equiv \\
\sum_{\substack{i\in I_0\\ \deg_{\rm coeff}\phi_i=k+\delta-e}}
c_i\tilde{\phi}_i^{(k+\delta-e)} \pmod{\pi}.
\end{multline}
Letting $\bar{\xi}^{(k)}\in\Omega^{m-1}_{{\bf F}_q[x]/{\bf F}_q}$
(resp.\ $\bar{\rho}_{i-1-j}^{(k+j-i)}\in\Omega^{m-1}_{{\bf
    F}_q[x]/{\bf F}_q}$) be the reduction mod $\pi$ of $\xi^{(k)}$
(resp.\ $\rho_{i-1-j}^{(k+j-i)}$), we see that equations (3.41) and
(3.42) imply
\begin{equation}
\sum_{j=0}^{i-1} df^{(\delta-j)}\wedge\bar{\rho}_{i-1-j}^{(k+j-i)} +
df^{(\delta-i)}\wedge\bar{\xi}^{(k)}=0 
\end{equation}
for $i=0,1,\ldots,e-2$ and
\begin{equation}
\sum_{j=0}^{e-2} df^{(\delta-j)}\wedge\bar{\rho}_{e-2-j}^{(k+j-e+1)} +
df^{(\delta-e+1)}\wedge\bar{\xi}^{(k)}=\sum_{\substack{i\in I_0\\ \deg_{\rm
      coeff}\phi_i=k+\delta-e}} \bar{c}_i\phi_i, 
\end{equation}
where $\bar{c}_i\in{\bf F}_q$ denotes the reduction mod $\pi$ of
$c_i$.  The definition of $E^{r,m-r}_e$ and the hypothesis that
$\{\phi_i\}_{\deg\phi_i=r}$ is linearly independent in $E^{r,m-r}_e$
then imply that $\bar{c}_i=0$ for all $i$.  This completes the proof
of Theorem 3.1. 

\section{Computation of $p$-adic cohomology}

We derive a series of corollaries to Theorem 3.1

\begin{corollary}
  In addition to the hypothesis of Theorem $3.1$, assume that the
  index set $I$ is finite, i.~e.,
\[ \dim_{{\bf F}_q}\biggl(\bigoplus_{r=0}^{\infty} E_e^{r,m-r}\biggr) <
  \infty. \]
Then the cohomology classes $\{[\tilde{\phi}_i]\}_{i\in I}$ span
  $H^m(\Omega^{\updot}_{C(b)},D)$ as $\tilde{\Omega}_0$-vector space.
\end{corollary}

{\it Proof}.  Let $\omega\in\Omega^m_{C(b)}$.  Without loss of
generality, we may assume $\omega\in\hat{F}^0\Omega^m_{C(b)}$.  By
Theorem 3.1 (see equation (3.6)), there exist
\[ \{c_i^{(0)}\}_{i\in I}\subseteq{\mathcal O}_{\tilde{\Omega}_0},
\quad \omega_0\in\hat{F}^{0}\Omega^m_{C(b)}, \quad
\xi_0\in\hat{F}^{-Mb(e-1)}\Omega^{m-1}_{C(b)} \]
such that
\begin{equation}
\omega=\pi\omega_0 +
\sum_{i\in I}c_i^{(0)}\tilde{\phi}_i + D(\xi_0).
\end{equation}
Suppose that for some $t\geq 0$ we have found
\[ \{c_i^{(t)}\}_{i\in I}\subseteq{\mathcal O}_{\tilde{\Omega}_0},
\quad \omega_t\in\hat{F}^{0}\Omega^m_{C(b)}, \quad
\xi_t\in\hat{F}^{-Mb(e-1)}\Omega^{m-1}_{C(b)} \]
such that
\begin{equation}
\omega=\pi^{t+1}\omega_t +
\sum_{i\in I}c_i^{(t)}\tilde{\phi}_i + D(\xi_t)
\end{equation}
and such that
\[ c_i^{(t)}-c_i^{(t-1)}\in\pi^t{\mathcal O}_{\tilde{\Omega}_0}, \quad
\xi_t-\xi_{t-1}\in\hat{F}^{-Mb(e-1)+t}\Omega^{m-1}_{C(b)}. \]
Applying Theorem 3.1 to $\omega_t$, we see that there exist
\[ \{\tilde{c}_i^{(t)}\}_{i\in I}\subseteq{\mathcal
  O}_{\tilde{\Omega}_0}, \quad
\omega_{t+1}\in\hat{F}^{0}\Omega^m_{C(b)}, \quad 
\xi'_t\in\hat{F}^{-Mb(e-1)}\Omega^{m-1}_{C(b)} \]
such that
\begin{equation}
\omega_t=\pi\omega_{t+1} +
\sum_{i\in I}\tilde{c}_i^{(t)}\tilde{\phi}_i + D(\xi'_t).
\end{equation}
Put $\xi_{t+1}=\xi_t+\pi^{t+1}\xi'_t$,
$c_i^{(t+1)}=c_i^{(t)}+\pi^{t+1}\tilde{c}_i^{(t)}$.  Substituting
(4.4) into (4.3) gives
\begin{equation}
\omega=\pi^{t+2}\omega_{t+1} +
\sum_{i\in I}c_i^{(t+1)}\tilde{\phi}_i + D(\xi_{t+1})
\end{equation}
with
\[ c_i^{(t+1)}-c_i^{(t)}\in\pi^{t+1}{\mathcal O}_{\tilde{\Omega}_0},
\quad \xi_{t+1}-\xi_t\in\hat{F}^{-Mb(e-1)+t+1}\Omega^{m-1}_{C(b)}. \]
It follows that $\{\xi_t\}_{t=0}^{\infty}$ converges to an element
$\xi\in\hat{F}^{-Mb(e-1)}\Omega^{m-1}_{C(b)}$ and
$\{c_i^{(t)}\}_{t=0}^{\infty}$ converges to an element
$c_i\in{\mathcal O}_{\tilde{\Omega}_0}$ for $i\in I$
satisfying
\[ \omega=\sum_{i\in I} c_i\tilde{\phi}_i+D(\xi). \]
This establishes the corollary.

\begin{corollary}
Suppose there exist $e,m$ such that $E^{r,m-r}_e=0$ for all $r\geq 0$.
Then for every rational number $b$ satisfying $(1.14)$,
we have $H^m(\Omega^{\updot}_{C(b)},D)=0$.  
\end{corollary}

{\it Proof}.  The hypothesis of Theorem 3.1 is satisfied with
$I=\emptyset$, so the conclusion follows immediately from Corollary 4.1.

To ensure that the cohomology classes $\{[\tilde{\phi}_i]\}_{i\in I}$
form a basis, we need an additional hypothesis.
\begin{corollary}
In addition to the hypothesis of Corollary $4.1$, suppose that for all
$r\in{\bf Z}$
\begin{equation}
\hat{E}^{r,m-r}_{Mb(e-1)+1}=\hat{E}^{r,m-r}_{Mb(e-1)+2}=\cdots .
\end{equation}
Then the set $\{[\tilde{\phi}_i]\}_{i\in I}$ is a basis for
$H^m(\Omega^{\updot}_{C(b)},D)$.  
\end{corollary}

{\it Proof}.  By Corollary 4.1, we only need to check linear
independence.  Suppose we had a relation
\begin{equation}
\sum_{i\in I}c_i\tilde{\phi_i}=D(\xi),
\end{equation}
where $\{c_i\}_{i\in I}\subseteq\tilde{\Omega}_0$ and
$\xi\in\Omega^{m-1}_{C(b)}$.  If the $c_i$ were not all zero, then
after multiplication by a suitable power of $\pi$ we may assume that
$c_i\in{\mathcal O}_{\tilde{\Omega}_0}$ for all $i$ and
$c_i\not\in\pi{\mathcal O}_{\tilde{\Omega}_0}$ for some $i$.  Thus the
left-hand side of (4.9) lies in $\hat{F}^0\Omega^m_{C(b)}$ but not 
in~$\hat{F}^1\Omega^m_{C(b)}$.  Since the filtration $\hat{F}^{\updot}$
is exhaustive, there exists $r\geq 0$ such that
\[ \xi\in\hat{F}^{-Mb(e-1)-r}\Omega^{m-1}_{C(b)}. \]
Equation (4.9) then says that $[\sum_{i\in I}c_i\tilde{\phi}_i]=0$ in
$\hat{E}^{0,m}_{Mb(e-1)+r+1}$.  By (4.8), we have
$[\sum_{i\in I}c_i\tilde{\phi}_i]=0$ in $\hat{E}^{0,m}_{Mb(e-1)+1}$.
Theorem 3.1 now implies $c_i\equiv 0\pmod{\pi}$ for all $i$, a
contradiction.  Thus $c_i=0$ for all $i$.

\begin{corollary}
In addition to the hypothesis of Corollary $4.1$, suppose that
\begin{equation}
E^{r,m-1-r}_e=E^{r,m+1-r}_e=0
\end{equation}
for all $r\geq 0$.  Then the set $\{[\tilde{\phi}_i]\}_{i\in I}$ is
a basis for $H^m(\Omega^{\updot}_{C(b)},D)$.  In particular,
\[ \dim_{\tilde{\Omega}_0} H^m(\Omega^{\updot}_{C(b)},D) = \dim_{{\bf
    F}_q} \biggl(\bigoplus_{r=0}^{\infty} E_e^{r,m-r}\biggr). \]
\end{corollary}

{\it Proof}.  Condition (4.11) and Theorem 3.1 imply that
\[ \hat{E}^{r,m-1-r}_{Mb(e-1)+1}=\hat{E}^{r,m+1-r}_{Mb(e-1)+1}=0 \]
for all $r\in{\bf Z}$.  General properties of spectral sequences then
imply that (4.8) holds, and we conclude by Corollary 4.7.

{\it Proof of Theorem $1.13$}.  The first assertion of Theorem 1.13
follows from Corollary~4.6.  Suppose that (1.11) holds.  The
discussion in the introduction shows that (1.11) implies
\[ \dim_{{\bf F}_q}\biggl(\bigoplus_{r=0}^{\infty} E_e^{r,n-r}\biggr)
= M_f. \]
Since (1.11) also implies condition (4.11) for $m=n$,
it follows from Corollary~4.10 that  
\[ \dim_{\tilde{\Omega}_0}H^n(\Omega^{\updot}_{C(b)},D)=M_f. \]

\section{Proof of Theorem 1.18}

Throughout this section, we assume the hypothesis of Theorem 1.18.  Put
\[ Z^k=\{\omega\in\Omega^k_{{\bf F}_q[x]/{\bf F}_q} \mid
df^{(\delta)}\wedge\omega=0\}. \]
We leave it to the reader to check from the definition of the spectral
sequence (1.6) that the conclusion of Theorem 1.18 is equivalent to
the following assertion. 
\begin{proposition}
Let $k<n$ and let $\omega\in Z^k$ be a homogeneous form such that
\begin{equation}
df^{(\delta')}\wedge\omega=df^{(\delta)}\wedge\xi
\end{equation}
for some homogeneous form $\xi\in\Omega^k_{{\bf F}_q[x]/{\bf F}_q}$.
Then there exist homogeneous forms $\eta_1\in Z^{k-1}$,
$\eta_2\in\Omega^{k-1}_{{\bf F}_q[x]/{\bf F}_q}$, such that
\begin{equation}
\omega= df^{(\delta')}\wedge\eta_1+df^{(\delta)}\wedge\eta_2.
\end{equation}
\end{proposition}

Before starting the proof, we observe that Kita\cite{KI} has,
in effect, characterized the elements of $Z^k$.  For notational
convenience, we define a $1$-form $\Theta\in \Omega^1_{{\bf
    F}_q[x]/{\bf F}_q}$ by
\[ \Theta = f_1\cdots f_r\sum_{i=1}^r a_i\frac{df_i}{f_i}\quad
\biggl(=f_1\cdots f_r\frac{df^{(\delta)}}{f^{(\delta)}}\biggr), \]
and for any $l$-tuple $1\leq i_1<\cdots<i_l\leq r$ we define
\[ \Omega_{i_1\cdots i_l}=\frac{df_{i_1}}{f_{i_1}}\wedge\cdots\wedge
\frac{df_{i_l}}{f_{i_l}}, \]
a rational $l$-form with logarithmic poles along the divisor
$f_1\cdots f_r=0$ in ${\bf A}^n$.  Note that $\Theta\wedge
\Omega_{i_1\cdots i_l}$ has polynomial coefficients, hence lies in
$\Omega^{l+1}_{{\bf F}_q[x]/{\bf F}_q}$.
\begin{proposition}
Let $k<n$ and let $\omega\in Z^k$ be homogeneous.  Then
\begin{equation}
\omega=\Theta\wedge\biggl(\sum_{l=0}^{k-1} \sum_{1\leq
  i_1<\cdots<i_l\leq r} \Omega_{i_1\cdots i_l}\wedge\alpha_{i_1\cdots
  i_l}\biggr), 
\end{equation}
for some homogeneous forms $\alpha_{i_1\cdots
    i_l}\in\Omega^{k-1-l}_{{\bf F}_q[x]/{\bf F}_q}$. 
\end{proposition}

{\it Proof}.  The Euler relation and the hypothesis $(p,\delta)=1$
imply that $f^{(\delta)}$ lies in the ideal generated by the $\partial
f^{(\delta)}/\partial x_i$.  The equation $df^{(\delta)}\wedge\omega=0$ then
implies, by a standard result (\cite[Theorem 16.4]{MA}, that there
exists $\alpha\in\Omega^{k-1}_{{\bf F}_q[x]/{\bf F}_q}$ such that 
\begin{equation}
f^{(\delta)}\omega=df^{(\delta)}\wedge\alpha, 
\end{equation}
or, equivalently,
\begin{equation}
f_1\cdots f_r\omega=\Theta\wedge\alpha.
\end{equation}
This implies that $a_if_1\cdots\hat{f}_i\cdots
f_r\,df_i\wedge\alpha\in(f_i)$ for $i=1,\ldots,r$.
But $(p,a_i)=1$, and our hypothesis implies $f_j,f_i$ form a regular
sequence for $j\neq i$, hence $df_i\wedge\alpha\in(f_i)$
for $i=1,\ldots,r$.  It then follows from \cite[Proposition
2.2.3]{KI} that
\begin{equation}
\alpha=f_1\cdots f_r\sum_{l=0}^{k-1} \sum_{1\leq
  i_1<\cdots<i_l\leq r} \Omega_{i_1\cdots i_l}\wedge\alpha_{i_1\cdots
  i_l}, 
\end{equation}
for some homogeneous forms $\alpha_{i_1\cdots i_l}\in
\Omega^{k-1-l}_{{\bf F}_q[x]/{\bf F}_q}$.  (Although Kita works over
${\bf C}$, his proof is valid for any field.)  Substituting (5.8) into
(5.7) gives the proposition. 

The following lemma is the main technical tool for the proof of
Proposition~5.1.
\begin{lemma}
Suppose that $\{i_1,\ldots,i_l\}$ is a nonempty subset of
$\{1,\ldots,r\}$ with either $l\geq 2$ or $l=1$ and $a_{i_1}>1$ and
that $\beta\in\Omega^m_{{\bf F}_q[x]/{\bf F}_q}$ 
is a homogeneous $m$-form with $m\leq n-l-1$ such that
\[ df^{(\delta')}\wedge df_{i_1}\wedge\cdots\wedge df_{i_l}\wedge\beta\equiv 0
\pmod{(f_{i_1},\ldots,f_{i_l})}. \]
Then there exist homogeneous $(m-1)$-forms $\beta_j$ for
$j=0,1,\ldots,l$ and $\beta'_j$ for $j=1,\ldots,l$ such that
\[ \beta=df^{(\delta')}\wedge\beta_0 + \sum_{j=1}^l df_{i_j}\wedge\beta_j +
\sum_{j=1}^l f_{i_j}\beta'_j. \]
\end{lemma}

{\it Proof}.  Consider the expansion of $df^{(\delta')}\wedge
df_{i_1}\wedge\cdots\wedge df_{i_l}$ relative to the basis for
$\Omega^{l+1}_{{\bf F}_q[x]/{\bf F}_q}$ consisting of the $(l+1)$-fold
exterior products of the $dx_i$, $i=1,\ldots,n$.  By the theorem of
\cite{SA}, it suffices to show that the coefficients in this expansion
generate an ideal $I$ of depth $n-l$ in the quotient ring ${\bf
  F}_q[x]/(f_{i_1},\ldots,f_{i_l})$, i.~e., the only maximal
ideal containing $I$ is the one generated by $x_1,\ldots,x_n$.
Suppose there were some other maximal ideal ${\bf m}$
containing $I$.  Then ${\bf m}$ would correspond to a point in ${\bf
  A}^n$, other than the origin, which is a common zero of
$f_{i_1},\ldots,f_{i_l}$ and at which there is a linear relation
between the differentials $df^{(\delta')},df_{i_1},\ldots,df_{i_l}$.  Since
$f_{i_1}=\cdots=f_{i_l}=0$ is a smooth complete intersection in ${\bf
  A}^n$ except at the origin, the differentials
$df_{i_1},\ldots,df_{i_l}$ are independent, so $df^{(\delta')}$ must
be a linear combination of $df_{i_1},\ldots,df_{i_l}$ at this point.
But then the Euler relations for the polynomials
$f^{(\delta')},f_{i_1},\ldots,f_{i_l}$, together with the hypothesis
that $(p,\delta')=1$, imply that $f^{(\delta')}$ 
also vanishes at this point, contradicting the hypothesis that
$f^{(\delta')}=f_{i_1}=\cdots=f_{i_l}=0$ defines a smooth complete
intersection in ${\bf A}^n$ except at the origin.  

Let $\omega\in Z^k$ and suppose that for some $s$, $1\leq s\leq r$, and
$b_s,\ldots,b_r$, $1\leq b_i\leq a_i$ for $i=s,\ldots,r$, we have
\begin{equation}
\omega = f_1^{a_1}\cdots f_{s-1}^{a_{s-1}}f_s^{b_s}\cdots
  f_r^{b_r}\sum_{i=1}^r 
  a_i\frac{df_i}{f_i}\wedge\biggl(\sum_{l=0}^{k-1} \sum_{1\leq
  i_1<\cdots<i_l\leq r}\Omega_{i_1\cdots i_l}\wedge\alpha_{i_1\cdots
  i_l}\biggr) 
\end{equation}
for some homogeneous forms $\alpha_{i_1\cdots i_l}\in
\Omega^{k-1-l}_{{\bf F}_q[x]/{\bf F}_q}$.  

\begin{lemma}
If $b_s<a_s$, then there exists $\eta\in Z^{k-1}$ such that
\begin{multline}
\omega-df^{(\delta')}\wedge\eta = \\
f_1^{a_1}\cdots
  f_{s-1}^{a_{s-1}}f_s^{b_s}\cdots f_r^{b_r}\sum_{i=1}^r 
  a_i\frac{df_i}{f_i}\wedge\biggl(\sum_{l=0}^{k-1} \sum_{1\leq
  i_1<\cdots<i_l\leq r}\Omega_{i_1\cdots i_l}\wedge\alpha'_{i_1\cdots
  i_l}\biggr), 
\end{multline}
where all $\alpha'_{i_1\cdots i_l}$ are divisible by $f_s$. 
\end{lemma}

Observe that Proposition 5.4 implies that every $\omega\in Z^k$ can be
written in the form (5.10) with $s=1$ and $b_1=\cdots=b_r=1$.  When
the conclusion of Lemma~5.11 holds, we can factor out $f_s$ from each
$\alpha'_{i_1\cdots i_l}$ and replace $b_s$ by $b_s+1$ in~(5.12).  
Induction on $b_s$ and $s$ then allows us to replace $b_s,\ldots,b_r$
by $a_s,\ldots,a_r$, respectively, in (5.12), giving the following.  
\begin{corollary}
If $\omega\in Z^k$ with $k<n$, then there exists $\eta\in Z^{k-1}$
such that
\[ \omega-df^{(\delta')}\wedge\eta= df^{(\delta)}\wedge\biggl(\sum_{l=0}^{k-1}
\sum_{1\leq i_1<\cdots<i_l\leq r}\Omega_{i_1\cdots
  i_l}\wedge\alpha_{i_1\cdots i_l}\biggr) \] 
for some homogeneous forms $\alpha_{i_1\cdots i_l}\in
  \Omega^{k-1-l}_{{\bf F}_q[x]/{\bf F}_q}$. 
\end{corollary}

{\it Proof of Lemma $5.11$}.  Put
\begin{equation}
\begin{split}
\tilde{\omega} &= \omega/(f_1^{a_1-1}\cdots
  f_{s-1}^{a_{s-1}-1}f_s^{b_s-1}\cdots f_r^{b_r-1}) \\
&= \Theta\wedge\biggl(\sum_{l=0}^{k-1} \sum_{1\leq
  i_1<\cdots<i_l\leq r}\Omega_{i_1\cdots i_l}\wedge\alpha_{i_1\cdots
  i_l}\biggr). 
\end{split}
\end{equation}
We show there exists $\tilde{\eta}\in Z^{k-1}$ such that
\begin{equation}
\tilde{\omega}-df^{(\delta')}\wedge\tilde{\eta} =
  \Theta\wedge\biggl(\sum_{l=0}^{k-1} \sum_{1\leq 
  i_1<\cdots<i_l\leq r}\Omega_{i_1\cdots i_l}\wedge\alpha'_{i_1\cdots
  i_l}\biggr), 
\end{equation}
where all $\alpha'_{i_1\cdots i_l}$ are divisible by $f_s$.  
Lemma 5.11 follows from (5.15) by multiplication by $f_1^{a_1-1}\cdots
  f_{s-1}^{a_{s-1}-1}f_s^{b_s-1}\cdots f_r^{b_r-1}$.

Equation (5.2) implies
\[ df^{(\delta')}\wedge\tilde{\omega}=f_s^{a_s-b_s}\cdots f_r^{a_r-b_r}
\Theta\wedge\xi, \] 
so if $b_s<a_s$ we have
\begin{equation}
df^{(\delta')}\wedge\tilde{\omega}\equiv 0 \pmod{(f_s)}.
\end{equation}
We prove the existence of $\tilde{\eta}$ by descending induction
on $l$.  Since $\alpha_{i_1\cdots i_l}=0$ for $l>k-1$, we may assume
that for some $m$, $0\leq m\leq k-1$, we have
$\alpha_{i_1\cdots i_l}$ divisible by $f_s$ for all $l\geq
m+1$.  We show that we can choose $\tilde{\eta}$ so that the
$\alpha'_{i_1\cdots i_l}$ are divisible by $f_s$ for all $l\geq m$.

Fix an $m$-tuple $1\leq i_1<\cdots<i_m\leq r$ with
$s\not\in\{i_1,\ldots,i_m\}$, say, $i_t<s<i_{t+1}$.  Expand the
right-hand side of (5.14) using the definitions of $\Theta$
and $\Omega_{i_1\cdots i_l}$.  The term containing
$df_{i_1}\wedge\cdots\wedge df_s\wedge\cdots\wedge df_{i_m}$ is
\begin{multline}
df_{i_1}\wedge\cdots\wedge df_s\wedge\cdots\wedge df_{i_m}\wedge 
f_1\cdots \hat{f}_{i_1}\cdots\hat{f}_s\cdots\hat{f}_{i_m}\cdots f_r 
\biggl((-1)^ta_s\alpha_{i_1\cdots i_m} + \\ \sum_{j=1}^t
(-1)^{j-1}a_{i_j}\alpha_{i_1\cdots\hat{\imath}_j\cdots s\cdots
  i_m} + \sum_{j=t+1}^m (-1)^j a_{i_j}\alpha_{i_1\cdots s\cdots
  \hat{\imath}_j\cdots i_m}\biggr).  
\end{multline}
Using the induction hypothesis that $\alpha_{i_1\cdots i_l}$ is
divisible by $f_s$ for all $l\geq m+1$, one sees that the
term containing any other product $df_{i_1}\wedge\cdots\wedge df_{i_l}$
(for $l\geq 0$) on the right-hand side of (5.14) lies in the ideal
$(f_s,f_{i_1},\ldots,f_{i_m})$.  Equation (5.16) then implies
\begin{multline}
df^{(\delta')}\wedge df_{i_1}\wedge\cdots\wedge df_s\wedge\cdots\wedge
df_{i_m}\wedge \\ 
f_1\cdots \hat{f}_{i_1}\cdots\hat{f}_s\cdots\hat{f}_{i_m}\cdots f_r 
\biggl((-1)^ta_s\alpha_{i_1\cdots i_m} +\sum_{j=1}^t
(-1)^{j-1}a_{i_j}\alpha_{i_1\cdots\hat{\imath}_j\cdots s\cdots
  i_m} + \\
\sum_{j=t+1}^m (-1)^j a_{i_j}\alpha_{i_1\cdots s\cdots
  \hat{\imath}_j\cdots i_m}\biggr) 
\equiv 0 \pmod{(f_s,f_{i_1},\ldots,f_{i_m})}.
\end{multline}

Since $\omega$ is a $k$-form with $k\leq n-1$, we may assume $m\leq
n-2$.  The smooth complete intersection hypothesis implies that for
$j\not\in\{s,i_1,\ldots,i_m\}$, $f_j,f_s,f_{i_1},\ldots,f_{i_m}$ form
a regular sequence, hence (5.18) implies
\begin{multline}
df^{(\delta')}\wedge df_{i_1}\wedge\cdots\wedge df_s\wedge\cdots\wedge
df_{i_m}\wedge \\
\biggl((-1)^ta_s\alpha_{i_1\cdots i_m} +\sum_{j=1}^t
(-1)^{j-1}a_{i_j}\alpha_{i_1\cdots\hat{\imath}_j\cdots s\cdots
  i_m} + \sum_{j=t+1}^m (-1)^j a_{i_j}\alpha_{i_1\cdots s\cdots
  \hat{\imath}_j\cdots i_m}\biggr) \\
\equiv 0 \pmod{(f_s,f_{i_1},\ldots,f_{i_m})}.
\end{multline}
It now follows by Lemma 5.9 (since $b_s<a_s$ implies $a_s>1$) that
there exist homogeneous forms 
$\gamma^{(j)}_{i_1\cdots i_m}$, $\delta^{(j)}_{i_1\cdots i_m}$, such
that
\begin{multline}
\alpha_{i_1\cdots i_m}=\frac{(-1)^t}{a_s}\biggl(\sum_{j=1}^t
(-1)^ja_{i_j}\alpha_{i_1\cdots\hat{\imath}_j\cdots s\cdots i_m}
+ \sum_{j=t+1}^m (-1)^{j-1}
a_{i_j}\alpha_{i_1\cdots s\cdots\hat{\imath}_j\cdots s\cdots
  i_m}\biggr) + \\
df^{(\delta')}\wedge\gamma^{(0)}_{i_1\cdots i_m} +
df_s\wedge\gamma^{(s)}_{i_1\cdots i_m} + f_s\delta^{(s)}_{i_1\cdots
  i_m}+\sum_{j=1}^m df_{i_j}\wedge\gamma^{(i_j)}_{i_1\cdots i_m} +
\sum_{j=1}^m f_{i_j}\delta^{(i_j)}_{i_1\cdots i_m}.
\end{multline}

Such a formula holds for every $m$-tuple $i_1,\ldots,i_m$ not containing
$s$.  Substitute these expressions into 
\begin{equation}
\Theta\wedge\biggl(\sum_{1\leq i_1<\cdots<i_m\leq r}
  \Omega_{i_1\cdots i_m}\wedge\alpha_{i_1\cdots i_m}\biggr)
\end{equation}
and expand.  After this substitution, only alphas indexed by $m$-tuples
containing $s$ remain.  We leave it to the reader to check that for
such an $m$-tuple, say,
\[ 1\leq j_1<\cdots<j_t<s<j_{t+1}<\cdots<j_{m-1}\leq r, \]
the contribution to (5.21) is
\[ \Theta\wedge\frac{(-1)^t\Theta}{a_sf_1\cdots
  f_r}\wedge \Omega_{j_1\cdots j_{m-1}}\wedge \alpha_{j_1\cdots
  s\cdots j_{m-1}} = 0. \]
Thus after substitution from (5.20), expression (5.21) simplifies to
\begin{multline}
\Theta\wedge\biggl(\sum_{\substack{1\leq i_1<\cdots<i_m\leq r\\
    s\not\in\{i_1,\ldots,i_m\}}} \Omega_{i_1\cdots i_m}\wedge \\
\biggl(df^{(\delta')}\wedge\gamma^{(0)}_{i_1\cdots i_m} +
df_s\wedge\gamma^{(s)}_{i_1\cdots i_m} 
+ f_s\delta^{(s)}_{i_1\cdots i_m}+ \sum_{j=1}^m
  f_{i_j}\delta^{(i_j)}_{i_1\cdots i_m}\biggr)\biggr). 
\end{multline}
We thus take
\begin{equation}
\tilde{\eta}=(-1)^{m+1}\Theta\wedge 
\biggl(\sum_{\substack{1\leq i_1<\cdots<i_m\leq r\\
    s\not\in\{i_1,\ldots,i_m\}}}\Omega_{i_1\cdots
  i_m}\wedge\gamma^{(0)}_{i_1\cdots i_m}\biggr) \in Z^{k-1},
\end{equation}
which (by (5.14) and (5.22)) gives
\begin{multline}
\tilde{\omega}-df^{(\delta')}\wedge\tilde{\eta}=\Theta\wedge
    \biggl( \sum_{\substack{l=0 \\ l\neq m}}^{k-1} \sum_{1\leq
    i_1<\cdots<i_l\leq r} \Omega_{i_1\cdots i_l}
\wedge\alpha_{i_1\cdots i_l} + \\
\sum_{\substack{1\leq i_1<\cdots<i_m\leq r\\ s\not\in\{i_1,\ldots,i_m\}}}
\Omega_{i_1\cdots i_m}\wedge\biggl( df_s\wedge\gamma^{(s)}_{i_1\cdots
    i_m} + f_s\delta^{(s)}_{i_1\cdots i_m} + \sum_{j=1}^m
f_{i_j}\delta^{(i_j)}_{i_1\cdots i_m}\biggr)\biggr).
\end{multline}
We rewrite this as
\begin{multline}
\tilde{\omega}-df^{(\delta')}\wedge\tilde{\eta}=\Theta\wedge
    \biggl( \sum_{\substack{l=0 \\ l\neq m}}^{k-1} \sum_{1\leq
    i_1<\cdots<i_l\leq r} \Omega_{i_1\cdots i_l}
\wedge\alpha_{i_1\cdots i_l} + \\
\sum_{\substack{1\leq i_1<\cdots<i_m\leq r\\ s\not\in\{i_1,\ldots,i_m\}}}
\Omega_{i_1\cdots i_m}\wedge f_s\delta^{(s)}_{i_1\cdots i_m} + 
\sum_{\substack{1\leq i_1<\cdots<i_m\leq r\\
    s\not\in\{i_1,\ldots,i_m\}}} \Omega_{i_1\cdots s\cdots
    i_m} \wedge (-1)^{m-t} f_s\gamma^{(s)}_{i_1\cdots i_m} + \\
\sum_{\substack{1\leq i_1<\cdots<i_m\leq r\\
    s\not\in\{i_1,\ldots,i_m\}}} \sum_{j=1}^m
    \Omega_{i_1\cdots\hat{\imath}_j\cdots i_m}\wedge
    (-1)^{m-j}df_{i_j}\wedge\delta^{(i_j)}_{i_1\cdots i_m}\biggr).
\end{multline}

For $l\geq m+2$, the coefficient of $\Omega_{i_1\cdots i_l}$ on the
right-hand side of (5.25) equals $\alpha_{i_1\cdots i_l}$, which is
divisible by $f_s$ by the induction hypothesis.  For $l=m+1$, the
coefficient of $\Omega_{i_1\cdots i_{m+1}}$ equals $\alpha_{i_1\cdots
  i_{m+1}}$ if $s\not\in\{i_1,\ldots,i_{m+1}\}$, while the coefficient
of $\Omega_{i_1\cdots s\cdots i_m}$ equals $\alpha_{i_1\cdots s\cdots
  i_m} + (-1)^{m-t} f_s\gamma^{(s)}_{i_1\cdots i_m}$.  In either case, it
is divisible by $f_s$ by the induction hypothesis.  For $l=m$, the
coefficient of $\Omega_{i_1\cdots i_m}$ equals $0$ if
$s\in\{i_1,\ldots,i_m\}$ and equals $f_s\delta^{(s)}_{i_1\cdots i_m}$
if $s\not\in\{i_1,\ldots,i_m\}$.  Thus the coefficient of
$\Omega_{i_1\cdots i_l}$ is divisible by $f_s$ for all $l\geq m$, and
by induction the proof of Lemma 5.11 is complete.  

{\it Proof of Proposition $5.1$}.  Suppose $\omega\in Z^k$ satisfies
(5.2).  By Corollary 5.13 we may assume that
\begin{equation}
\omega=df^{(\delta)}\wedge\sum_{l=1}^{k-1}\sum_{1\leq i_1<\cdots<i_l\leq r}
\Omega_{i_1\cdots i_l}\wedge\alpha_{i_1\cdots i_l}.
\end{equation}
Note that we may start the outer sum at $l=1$ rather than at $l=0$.  
We prove by induction on $s$ that for $1\leq s\leq r$ we can find
$\eta_1\in Z^{k-1}$, $\eta_2\in\Omega^{k-1}_{{\bf F}_q[x]/{\bf F}_q}$,
such that
\begin{equation}
\omega-df^{(\delta')}\wedge\eta_1-df^{(\delta)}\wedge\eta_2 = df^{(\delta)}
\wedge\sum_{l=1}^{k-1}\sum_{1\leq i_1<\cdots<i_l\leq r}
\Omega_{i_1\cdots i_l}\wedge\alpha'_{i_1\cdots i_l}
\end{equation}
with all $\alpha'_{i_1\cdots i_l}$ divisible by $f_1\cdots f_s$.
When $s=r$, $\Omega_{i_1\cdots i_l}\wedge\alpha'_{i_1\cdots i_l}$ has
polynomial coefficients for all $i_1,\ldots,i_l$, so equation (5.27)
establishes Proposition 5.1.

So suppose that for some $s$, $1\leq s\leq r$, all $\alpha_{i_1\cdots
  i_l}$ in (5.26) are divisible by $f_1\cdots f_{s-1}$.  In this case, if
$i_l<s$, then $\Omega_{i_1\cdots i_l}\wedge\alpha_{i_1\cdots
  i_l}\in\Omega^{k-1}_{{\bf F}_q[x]/{\bf F}_q}$, so by replacing
$\omega$ by $\omega-df^{(\delta)}\wedge\Omega_{i_1\cdots
  i_l}\wedge\alpha_{i_1\cdots i_l}$, we may assume
\begin{equation}
\alpha_{i_1\cdots i_l}=0 \qquad\text{when $i_l<s$.}
\end{equation}
We prove we can choose
$\eta_1,\eta_2$ so that all $\alpha'_{i_1\cdots i_l}$ are divisible by
$f_1\cdots f_s$ by descending induction on $l$.  Since
$\alpha_{i_1\cdots i_l}=0$ for $l>k-1$, we may assume that for some
$m$, $1\leq m\leq k-1$, $\alpha_{i_1\cdots i_l}$ is divisible by
$f_1\cdots f_s$ for $l\geq m+1$.  We show that we can choose
$\eta_1,\eta_2$ so that $\alpha'_{i_1\cdots i_l}$ is divisible by
$f_1\cdots f_s$ for $l\geq m$.

Put
\begin{equation}
\begin{split}
\tilde{\omega} &= \omega/(f_1^{a_1-1}\cdots f_r^{a_r-1}) \\
&= \Theta\wedge\biggl(\sum_{l=1}^{k-1} \sum_{\substack{1\leq
  i_1<\cdots<i_l\leq r\\ s\leq i_l}}\Omega_{i_1\cdots
i_l}\wedge\alpha_{i_1\cdots i_l}\biggr). 
\end{split}
\end{equation}
We show there exist $\tilde{\eta}_1\in Z^{k-1}$,
$\tilde{\eta}_2\in\Omega^{k-1}_{{\bf F}_q[x]/{\bf F}_q}$ such that
\begin{equation}
\tilde{\omega}-df^{(\delta')}\wedge\tilde{\eta}_1-
\Theta\wedge\tilde{\eta}_2 = \Theta
\wedge\sum_{l=1}^{k-1}\sum_{1\leq i_1<\cdots<i_l\leq r}
\Omega_{i_1\cdots i_l}\wedge\alpha'_{i_1\cdots i_l} 
\end{equation}
with all $\alpha'_{i_1\cdots i_l}$ divisible by $f_1\cdots f_s$.
Equation (5.27) follows from (5.30) by multiplication by $f_1^{a_1-1}\cdots
f_r^{a_r-1}$.  

Equation (5.2) implies
\begin{equation}
df^{(\delta')}\wedge\tilde{\omega}\equiv 0 \pmod{(f_s,f_i)}
\end{equation}
for all $i\neq s$.  Fix an $m$-tuple $1\leq i_1<\cdots<i_m\leq r$ with
$s\not\in\{i_1,\ldots,i_m\}$ and $s<i_m$, say, $i_t<s<i_{t+1}$.  The
term containing $df_{i_1}\wedge\cdots\wedge df_s\wedge\cdots\wedge
df_{i_m}$ in the expansion of the right-hand side of (5.30) is given by
(5.17).  Using the induction hypothesis that $\alpha_{i_1\cdots i_l}$ is
divisible by $f_1\cdots f_s$ for all $l\geq m+1$ and by $f_1\cdots
f_{s-1}$ for all $l\geq 1$, one sees that the
term containing any other product $df_{i_1}\wedge\cdots\wedge df_{i_l}$
(for $l\geq 1$) on the right-hand side of (5.30) lies in the ideal
$(f_s,\{f_{i_j}\}_{i_j>s},\{f_{i_j}^2\}_{i_j<s})$.  Since $\{i_j\mid
i_j>s\}\neq\emptyset$, equation (5.31) implies that
\begin{multline}
df^{(\delta')}\wedge df_{i_1}\wedge\cdots\wedge df_s\wedge\cdots\wedge
df_{i_m}\wedge \\ 
f_1\cdots \hat{f}_{i_1}\cdots\hat{f}_s\cdots\hat{f}_{i_m}\cdots f_r 
\biggl((-1)^ta_s\alpha_{i_1\cdots i_m} +\sum_{j=1}^t
(-1)^{j-1}a_{i_j}\alpha_{i_1\cdots\hat{\imath}_j\cdots s\cdots
  i_m} + \\
\sum_{j=t+1}^m (-1)^j a_{i_j}\alpha_{i_1\cdots s\cdots
  \hat{\imath}_j\cdots i_m}\biggr) 
\equiv 0 \pmod{(f_s,\{f_{i_j}\}_{i_j>s},\{f_{i_j}^2\}_{i_j<s})}.
\end{multline}
As noted earlier, we may assume $m\leq n-2$, which implies that
$f_j,f_s,f_{i_1},\ldots,f_{i_m}$ is a regular sequence.  We then
deduce from (5.32) that
\begin{multline}
df^{(\delta')}\wedge df_{i_1}\wedge\cdots\wedge df_s\wedge\cdots\wedge
df_{i_m}\wedge \\
\biggl((-1)^ta_s\frac{\alpha_{i_1\cdots i_m}}{f_1\cdots f_{s-1}}+ 
\sum_{j=1}^t (-1)^{j-1}a_{i_j}
  \frac{\alpha_{i_1\cdots\hat{\imath}_j\cdots s\cdots i_m}}{f_1\cdots
    f_{s-1}} + \sum_{j=t+1}^m (-1)^j a_{i_j}\frac{\alpha_{i_1\cdots s\cdots
  \hat{\imath}_j\cdots i_m}}{f_1\cdots f_{s-1}} \biggr) \\
\equiv 0 \pmod{(f_s,f_{i_1},\ldots,f_{i_m})}.
\end{multline}
Applying Lemma 5.9 (which is permissible since $m\geq 1$) and then
multiplying by $f_1\cdots f_{s-1}$, we see 
that there exist homogeneous forms
$\gamma^{(j)}_{i_1\cdots i_m}$, $\delta^{(j)}_{i_1\cdots i_m}$, all
divisible by $f_1\cdots f_{s-1}$, such that
\begin{multline}
\alpha_{i_1\cdots i_m}=\frac{(-1)^t}{a_s}\biggl(\sum_{j=1}^t
(-1)^ja_{i_j}\alpha_{i_1\cdots\hat{\imath}_j\cdots s\cdots i_m}
+ \sum_{j=t+1}^m (-1)^{j-1}
a_{i_j}\alpha_{i_1\cdots s\cdots\hat{\imath}_j\cdots s\cdots
  i_m}\biggr) + \\
df^{(\delta')}\wedge\gamma^{(0)}_{i_1\cdots i_m} +
df_s\wedge\gamma^{(s)}_{i_1\cdots i_m} + f_s\delta^{(s)}_{i_1\cdots
  i_m}+\sum_{j=1}^m df_{i_j}\wedge\gamma^{(i_j)}_{i_1\cdots i_m} +
\sum_{j=1}^m f_{i_j}\delta^{(i_j)}_{i_1\cdots i_m}.
\end{multline}

Such a formula holds for every $m$-tuple $i_1,\ldots,i_m$ not containing
$s$ with $s<i_m$.  Substitute these expressions into 
\begin{equation}
\Theta\wedge\biggl(\sum_{\substack{1\leq i_1<\cdots<i_m\leq
  r\\ s\leq i_m}}\Omega_{i_1\cdots i_m}\wedge\alpha_{i_1\cdots i_m}\biggr).
\end{equation}
After this substitution, only alphas indexed by $m$-tuples
containing $s$ remain.  Consider such an $m$-tuple, say,
\[ 1\leq j_1<\cdots<j_t<s<j_{t+1}<\cdots<j_{m-1}\leq r. \]
If $t<m-1$, the reader may check that the contribution to (5.35) is
\[ \Theta\wedge\frac{(-1)^t\Theta}{a_sf_1\cdots
  f_r} \wedge \Omega_{j_1\cdots j_{m-1}}\wedge \alpha_{j_1\cdots
  s\cdots j_{m-1}} = 0. \]
In the case $t=m-1$ (i.~e., $j_i<s$ for all $i$), the reader may check
  that the contribution is 
\begin{multline*}
\Theta\wedge\frac{(-1)^{m-1}}{a_s}\sum_{i=s}^r
a_i\Omega_{ij_1\cdots j_{m-1}}\wedge\alpha_{j_1\cdots j_{m-1}s} \\
= \Theta\wedge\frac{(-1)^m}{a_s}
\sum_{i=1}^{s-1} a_i\Omega_{ij_1\cdots j_{m-1}}\wedge \alpha_{j_1\cdots
  j_{m-1}s},
\end{multline*}
since $\sum_{i=1}^r a_i\Omega_{ij_1\cdots j_{m-1}}=\Theta
\wedge\Omega_{j_1\cdots j_{m-1}}$.  Put
\[ \tilde{\eta}_2= \frac{(-1)^m}{a_s}
\sum_{i=1}^{s-1} a_i\Omega_{ij_1\cdots j_{m-1}}\wedge \alpha_{j_1\cdots
  j_{m-1}s}. \]
For $i=1,\ldots,s-1$, we have $\{i,j_1,\ldots,j_{m-1}\}\subseteq
\{1,\ldots,s-1\}$, and since $\alpha_{j_1\cdots j_{m-1}s}$ is
divisible by $f_1\cdots f_{s-1}$, we conclude that $\tilde{\eta}_2\in
\Omega^{k-1}_{{\bf F}_q[x]/{\bf F}_q}$.  

Thus after substitution from (5.34), expression (5.35) simplifies to
\begin{multline}
\Theta\wedge 
\biggl(\tilde{\eta}_2 + \sum_{\substack{1\leq i_1<\cdots<i_m\leq r\\
    s\not\in\{i_1,\ldots,i_m\}\\ s<i_m}}\Omega_{i_1\cdots i_m}\wedge \\
\biggl(df^{(\delta')}\wedge\gamma^{(0)}_{i_1\cdots i_m} +
df_s\wedge\gamma^{(s)}_{i_1\cdots i_m} 
+ f_s\delta^{(s)}_{i_1\cdots i_m}+ \sum_{j=1}^m
  f_{i_j}\delta^{(i_j)}_{i_1\cdots i_m}\biggr)\biggr). 
\end{multline}
We now take
\begin{equation}
\tilde{\eta}_1=(-1)^{m+1}\Theta\wedge 
\biggl(\sum_{\substack{1\leq i_1<\cdots<i_m\leq r\\
    s\not\in\{i_1,\ldots,i_m\}\\ s<i_m}}\Omega_{i_1\cdots
  i_m}\wedge\gamma^{(0)}_{i_1\cdots i_m}\biggr) \in Z^{k-1}
\end{equation}
and conclude (see (5.29))
\begin{multline}
\tilde{\omega}-df^{(\delta')}\wedge\tilde{\eta}_1 -
\Theta\wedge\tilde{\eta}_2 =  
\Theta\wedge \biggl( \sum_{\substack{l=0 \\ l\neq
    m}}^{k-1} \sum_{\substack{1\leq i_1<\cdots<i_l\leq r\\ s\leq i_l}}
\Omega_{i_1\cdots i_l}\wedge\alpha_{i_1\cdots i_l} + \\
\sum_{\substack{1\leq i_1<\cdots<i_m\leq r\\
    s\not\in\{i_1,\ldots,i_m\} \\ s<i_m}}\Omega_{i_1\cdots i_m}
\wedge\biggl( df_s\wedge\gamma^{(s)}_{i_1\cdots i_m} +
f_s\delta^{(s)}_{i_1\cdots i_m} + \sum_{j=1}^m
f_{i_j}\delta^{(i_j)}_{i_1\cdots i_m}\biggr)\biggr).
\end{multline}
Rewriting the right-hand side of (5.38) as in (5.25) and arguing as we
did in the proof of Lemma 5.11 shows that the coefficient of
$\Omega_{i_1\cdots i_l}$ is divisible by $f_1\cdots f_s$ for $l\geq
m$.  This completes the proof of Proposition 5.1.

\end{document}